\documentclass[11pt]{article}
\pdfoutput=1
\usepackage{hyperref}
\usepackage{url}
\usepackage{amsmath}
\usepackage{bm}
\usepackage{amsfonts}
\usepackage{algorithm}
\usepackage{algorithmic}
\usepackage{amsthm,paralist}
\usepackage{enumitem,graphicx}
\usepackage{color}
\usepackage{epstopdf}
\usepackage{wrapfig}
\usepackage{subcaption} 
\usepackage{amssymb}

\DeclareMathOperator{\argmin}{\mbox{argmin}}
\def\bn{\hfill \\ \smallskip\noindent}
\def\argmin{\mathop{\rm argmin}}

\def\prox{\mbox{prox}}

\def\dist{\mbox{dist\,}}

\newcommand{\beq}{\begin{equation}}
\newcommand{\eeq}{\end{equation}}
\newcommand{\st}{{\rm s.t.}}

\begin{document}
\def\pn {\par\smallskip\noindent}
\def \bn {\hfill \\ \smallskip\noindent}
\newcommand{\fs}{f_1,\ldots,f_s}
\newcommand{\f}{\vec{f}}
\newcommand{\hf}{\hat{f}}
\newcommand{\hx}{\hat{x}}
\newcommand{\hy}{\hat{y}}
\newcommand{\hz}{\hat{z}}
\newcommand{\hw}{\hat{w}}
\newcommand{\tw}{\tilde{w}}
\newcommand{\hlambda}{\hat{\lambda}}
\newcommand{\hbeta}{\hat{\beta}}
\newcommand{\tG}{\widetilde{G}}
\newcommand{\tg}{\widetilde{g}}
\newcommand{\barhx}{\bar{\hat{x}}}
\newcommand{\vecx}{x_1,\ldots,x_m}
\newcommand{\xoy}{x\rightarrow y}
\newcommand{\barx}{{\bar x}}
\newcommand{\bary}{{\bar y}}
\newcommand{\hrho}{\widehat{\rho}}
\newtheorem{theorem}{Theorem}[section]
\newtheorem{lemma}{Lemma}[section]
\newtheorem{corollary}{Corollary}[section]
\newtheorem{proposition}{Proposition}[section]
\newtheorem{definition}{Definition}[section]
\newtheorem{claim}{Claim}[section]
\newtheorem{remark}{Remark}[section]

\def\smskip{\par\vskip 5 pt}
\def\proof{\bn {\bf Proof.} }
\def\QED{\hfill{\bf Q.E.D.}\smskip}
\def\qed{\quad{\bf q.e.d.}\smskip}
\newcommand{\cE}{\mathcal{E}}
\newcommand{\cM}{\mathcal{M}}
\newcommand{\cN}{\mathcal{N}}
\newcommand{\cJ}{\mathcal{J}}
\newcommand{\cT}{\mathcal{T}}
\newcommand{\bx}{\mathbf{x}}
\newcommand{\bp}{\mathbf{p}}
\newcommand{\bX}{\mathbf{X}}
\newcommand{\bY}{\mathbf{Y}}
\newcommand{\bP}{\mathbf{P}}
\newcommand{\bD}{\mathbf{D}}
\newcommand{\bA}{\mathbf{A}}
\newcommand{\bB}{\mathbf{B}}
\newcommand{\bfM}{\mathbf{M}}
\newcommand{\bL}{\mathbf{L}}
\newcommand{\bz}{\mathbf{z}}
\newcommand{\cF}{\mathcal{F}}
\newcommand{\bF}{\mathbf{F}}
\newcommand{\cR}{\mathcal{R}}
\newcommand{\bzero}{\mathbf{0}}

\newcommand{\blue}{\color{blue}}
\newcommand{\red}{\color{red}}

\title{NESTT: A Nonconvex Primal-Dual Splitting Method for Distributed and Stochastic Optimization}

\author{{Davood Hajinezhad, Mingyi Hong }\thanks{ Department of Industrial \& Manufacturing Systems Engineering and Department of Electrical \& Computer Engineering, Iowa State University,  Emails: \texttt{\{dhaji,mingyi\}@iastate.edu}} \and {Tuo Zhao}\thanks{Department of Computer Science,  Johns Hopkins University, Email: \texttt{tzhao5@jhu.edu}} \and{Zhaoran Wang}\thanks{ Department of Operations Research and Financial Engineering,  Princeton University, Email: \texttt{zhaoran@princeton.edu} 
}}
\date{}

	\maketitle
	\vspace{-0.5cm}
	\begin{abstract}
		We study a stochastic and distributed algorithm for nonconvex  problems whose objective consists of a sum of $N$ nonconvex $L_i/N$-smooth functions, plus a  nonsmooth regularizer. The proposed NonconvEx primal-dual SpliTTing (NESTT) algorithm splits the problem into $N$ subproblems, and utilizes an augmented Lagrangian based primal-dual scheme to solve it in a distributed and stochastic manner. With a special non-uniform sampling, a version of NESTT achieves $\epsilon$-stationary solution  using $\mathcal{O}((\sum_{i=1}^N\sqrt{L_i/N})^2/\epsilon)$ gradient evaluations, which can be up to $\mathcal{O}(N)$ times better than the (proximal) gradient descent methods. It also achieves Q-linear convergence rate for nonconvex $\ell_1$ penalized quadratic problems with polyhedral constraints. Further, we reveal  a fundamental connection between {\it primal-dual} based methods and a few {\it primal only} methods such as IAG/SAG/SAGA. 
	\end{abstract}

\section{Introduction}
Consider the following nonconvex and nonsmooth constrained optimization problem
\begin{align}\label{eq:original}
\min_{z\in Z}\quad  f(z):=  \frac{1}{N}\sum_{i=1}^{N} g_i(z) +g_0(z)+ p(z),
\end{align}
where $Z\subseteq \mathbb{R}^d$; for each $i\in\{0,\cdots, N\}$, $g_i:\mathbb{R}^d \to \mathbb{R}$ is a smooth possibly nonconvex function which has $L_i$-Lipschitz continuous gradient; $p(z):\mathbb{R}^d\to\mathbb{R}$ is a lower semi-continuous convex but possibly nonsmooth function. Define $g(z): =\frac{1}{N} \sum_{i=1}^{N} g_i(z)$ for notational simplicity.

Problem \eqref{eq:original} is quite general. It arises frequently in applications such as machine learning and signal processing; see a recent survey \cite{Cevher14}. In particular, each smooth functions $\{g_i\}_{i=1}^{N}$ can represent: 1) a mini-batch of loss functions modeling data fidelity, such as the $\ell_2$ loss, the logistic loss, etc; 2) nonconvex activation functions for neural networks, such as the $\rm{logit}$  or the $\tanh$ functions; 3) nonconvex utility functions used in signal processing, machine learning, and resource allocation, see \cite{bjornson13}, and \cite{Hajinezhad16}. The smooth function $g_0$ can represent smooth nonconvex regularizers such as the non-quadratic penalties \cite{Antoniadis2009}, or the smooth part of the SCAD or MCP regularizers (which is a concave function) \cite{wang14nonconvex}. The convex function $p$ can take the following form: 1) nonsmooth convex regularizers such as $\ell_1$ and $\ell_2$ functions; 2) an indicator function for convex and closed feasible set $Z$, denoted as $\iota_{Z}(\cdot)$; 3) convex functions without global Lipschitz continuous gradient, such as $p(z)=z^4$ or $p(z)=1/z + \iota_{z\ge 0}(z)$.

In this work we solve \eqref{eq:original} in a stochastic and distributed manner. We consider the setting in which $N$ distributed agents each having the knowledge of one smooth function $\{g_i\}_{i=1}^{N}$, and they are connected to a cluster center which handles $g_0$ and $p$. At any given time, a randomly selected agent is activated and performs computation to optimize its local objective. Such distributed computation model has been popular in large-scale machine learning and signal processing \cite{BoydADMMsurvey2011}. Such model is also closely related to the (centralized) stochastic {\it finite-sum} optimization problem \cite{Lan15,Defazio14,Johnson13,reddi16fast,zhu16VR,Schmidt13}, in which each time the iterate is updated based on the gradient information of a random component function. One of the key differences between these two problem types is that in the distributed setting there can be disagreement between local copies of the optimization variable $z$, while in the centralized setting only one copy of $z$ is maintained.

{\bf Our Contributions.} We propose a class of NonconvEx primal-dual SpliTTing (NESTT) algorithms for problem \eqref{eq:original}. We split $z\in\mathbb{R}^d$ into local copies of $x_i\in\mathbb{R}^d$, while enforcing the equality constraints $x_i = z$ for all $i$. That is, we consider the following reformulation of \eqref{eq:original}
\begin{align}\label{eq:distributed}
\min_{x,z\in\mathbb{R}^d}\quad  \ell(x,z):=  \frac{1}{N}\sum_{i=1}^{N} g_i(x_i) + g_0(z) + h(z), \quad \st\; x_i=z, \; i=1,\cdots, N,
\end{align}
where $h(z):=\iota_{Z}(z) + p(z)$, $x:=[x_1; \cdots ;x_N]$. 
Our algorithm uses the Lagrangian relaxation of the equality constraints, and at each iteration a (possibly non-uniformly) randomly selected primal variable is optimized, followed by an approximate dual ascent step. Note that such splitting scheme has been popular in the convex setting \cite{BoydADMMsurvey2011}, but not so when the problem becomes nonconvex.  

The NESTT is one of the first stochastic algorithms for distributed nonconvex nonsmooth optimization, with provable and nontrivial convergence rates. Our main contribution is given below. First, in terms of some primal and dual optimality gaps, NESTT converges sublinearly to a point belongs to stationary solution set of \eqref{eq:distributed}. Second, NESTT converges Q-linearly for certain nonconvex $\ell_1$ penalized quadratic problems. To the best of our knowledge, this is the first time that linear convergence is established for stochastic and distributed optimization of such type of problems. Third, we show that a gradient-based NESTT with {\it non-uniform sampling}  achieves an $\epsilon$-stationary solution of \eqref{eq:original} using $\mathcal{O}((\sum_{i=1}^{N}\sqrt{L_i/N})^2 /\epsilon)$ gradient evaluations.
Compared with the classical gradient descent, which in the worst case requires $\mathcal{O}(\sum_{i=1}^{N} L_i/\epsilon)$ gradient evaluation to achieve $\epsilon$-stationarity \cite{Nesterov04}, our obtained rate can be up to $\mathcal{O}(N)$ times better in the case where the  $L_i$'s are not equal.

Our work also reveals a fundamental connection between {\it primal-dual} based algorithms and the {\it primal only} average-gradient based algorithm such as SAGA/SAG/IAG \cite{Defazio14,reddi16fast,Schmidt13,Blatt07}. With the key observation that the dual variables in NESTT serve as the ``memory" of the past gradients, one can specialize NESTT to SAGA/SAG/IAG. Therefore, NESTT naturally generalizes these algorithms to the nonconvex nonsmooth setting. It is our hope that by bridging the primal-dual splitting algorithms and primal-only algorithms (in {\it both} the convex and nonconvex setting), there can be significant further research developments benefiting both algorithm classes.

\noindent{\bf Related Work.} Many  stochastic algorithms have been designed for \eqref{eq:distributed} when it is convex. In these algorithms the component functions $g_i$'s are randomly sampled and optimized.  Popular algorithms include the SAG/SAGA \cite{Defazio14,Schmidt13}, the SDCA \cite{shalev14proximaldual}, the SVRG \cite{Johnson13}, the RPDG \cite{Lan15}  and so on. When the problem becomes nonconvex, the well-known incremental based algorithm can be used \cite{survrit12,bertsekas2000incremental}, but these methods generally lack convergence rate guarantees. The SGD based method has been studied in \cite{Ghadimi13:stochastic_nonconvex}, with $\mathcal{O}(1/\epsilon^2)$ convergence rate. Recent works \cite{zhu16VR} and \cite{reddi16fast} develop algorithms based on SVRG and SAGA for a special case of \eqref{eq:original} where the entire problem is smooth and unconstrained. To the best of our knowledge there has been no stochastic algorithms with provable, and non-trivial, convergence rate guarantees for solving problem \eqref{eq:original}. On the other hand, distributed stochastic algorithms for solving problem \eqref{eq:original} in the nonconvex setting has been proposed in \cite{hong14nonconvex_admm}, \cite{Hajinezhad15}, in which each time a randomly picked subset of agents update their local variables. However there has been no convergence rate analysis for such distributed stochastic scheme. There has been some recent distributed algorithms designed for \eqref{eq:original} \cite{Lorenzo16}, but again without global convergence rate guarantee.

\noindent{\bf Preliminaries.} The augmented Lagrangian function for problem \eqref{eq:original} is given by:
{
	\begin{align}\label{eq:aug}
	L\left(x, z; \lambda\right) =\sum_{i=1}^{N}\left(\frac{1}{N}g_i(x_i)+ \langle \lambda_i, x_i-z\rangle+\frac{\eta_i}{2}\|x_i-z\|^2\right) + g_0(z) + h(z),
	\end{align}}
where $\lambda:=\{\lambda_i\}_{i=1}^{N}$ is the set of dual variables, and  $\eta:=\{\eta_i>0\}_{i=1}^{N}$ are penalty parameters.

We make the following assumptions about problem \eqref{eq:original} and the function \eqref{eq:aug}.
\begin{enumerate}
	\item [A-(a)] The function $f(z)$ is bounded from below over $Z\cap \hbox{int}(\hbox{dom } f)$:
	$\underbar{f} : =\min_{z\in Z} f(z) >-\infty.$ $p(z)$ is a convex lower semi-continuous function; $Z$ is a closed convex set.
	\item [A-(b)] The $g_i$'s and $g$ have Lipschitz continuous gradients, i.e.,
	\begin{align*}
	\|\nabla g(y)-\nabla g(z)\|\le L\|y-z\|, \text{and} \quad \|\nabla g_i(y)-\nabla g_i(z)\|\le L_i\|y-z\|, \; \forall~y,z
	\end{align*}
	Clearly $L\le 1/N\sum_{i=1}^{N}L_i$, and the equality can be achieved in the worst case. For simplicity of analysis we will further assume that $L_0\le \frac{1}{N}\sum_{i=1}^{N} L_i.$
	\item [A-(c)] Each $\eta_i$ in \eqref{eq:aug} satisfies $\eta_i > L_i/N$; if $g_0$ is nonconvex, then $\sum_{i=1}^{N}\eta_i > 3 L_0$.
\end{enumerate}
Assumption A-(c) implies that $L\left(x, z; \lambda\right)$ is {\it strongly convex} w.r.t. each $x_i$ and $z$, with modulus  $\gamma_i := \eta_i - L_i/N$ and $\gamma_{z} = \sum_{i=1}^{N}\eta_i- L_0$, respectively \cite[Theorem 2.1]{ZLOBEC05}.

We then define the {\it prox-gradient} ({\texttt{pGRAD}}) for  \eqref{eq:original}, which will serve as a measure of stationarity. It can be checked that the {\texttt{pGRAD}} vanishes at the set of stationary solutions of  \eqref{eq:original} \cite{meisam14nips}.
\begin{definition}
	The \textit{proximal gradient} of problem \eqref{eq:original} is given by (for any $\gamma>0$){\small
		\[
		\tilde \nabla f_{\gamma}(z):=\gamma\left(z-{\rm \prox}^{\gamma}_{p+\iota_{Z}}[z-1/\gamma\nabla (g(z)+g_0(z))]\right), \; \mbox{with} \quad  {\rm \prox}^{\gamma}_{p+\iota_{Z}}[u]:=\argmin_{u\in Z}\;\;p(u)+\frac{\gamma}{2}\|z-u\|^2.
		\]}
\end{definition}
\section{The NESTT-G Algorithm}
\noindent{\bf Algorithm Description.} We present a primal-dual splitting scheme for the reformulated problem  \eqref{eq:distributed}. The algorithm is referred to as the NESTT with Gradient step (NESTT-G) since each agent only requires to know the gradient of each component function. To proceed, let us define the following function (for some constants $\{\alpha_i>0\}_{i=1}^{N}$):
{
	\begin{align*}
	V_i(x_i,z; \lambda_i) = {\frac{1}{N}} g_i(z) + {\frac{1}{N}}\langle\nabla g_{i}(z), x_{i}-z\rangle+\langle \lambda_{i}, x_{i}-z\rangle+\frac{\alpha_i\eta_i}{2}\|x_{i}-z\|^2.
	\end{align*}}
Note that $V_i(\cdot)$ is related to $L(\cdot)$ in the following way: it is a quadratic approximation (approximated at the point $z$) of $L(x,y;\lambda)$ w.r.t. $x_i$. The parameters $\alpha:= \{\alpha_i\}_{i=1}^{N}$ give some freedom to the algorithm design, and they are critical in  improving convergence rates as well as in establishing connection between NESTT-G with a few primal only stochastic optimization schemes. 

The algorithm proceeds as follows.  Before each iteration begins the cluster center broadcasts $z$  to everyone. At iteration $r+1$ a randomly selected agent $i_{r}\in\{1,2,\cdots N\}$ is picked, who  minimizes $V_{i_{r}}(\cdot)$ w.r.t. its local variable $x_{i_{r}}$, followed by a dual ascent step for $\lambda_{i_{r}}$.  The rest of the agents update their local variables by simply setting them to $z$. The cluster center then minimizes $L(x, z; \lambda)$ with respect to $z$.  See Algorithm 1 for details.
\begin{algorithm}[tb]\label{alg:nesttg}
	\caption{NESTT-G Algorithm}
	\begin{algorithmic}[1]
		\FOR{$r=1$ {\bfseries to} $R$}
		\STATE Pick $i_{r}\in \{1,2,\cdots, N\}$ with probability $p_{i_{r}}$ and update $(x,\lambda)$
		\begin{flalign}
		x^{r+1}_{i_{r}}&=\arg\min_{x_{i_{r}}} V_{i_{r}}\left(x_{i_{r}},z^{r}, \lambda_{i_{r}}^r\right);&&\label{eq:x_i:nestt1}\\
		\lambda^{r+1}_{i_{r}}&=\lambda_{i_{r}}^{r}+\alpha_{i_{r}}\eta_{i_{r}}\left(x^{r+1}_{i_{r}}-z^{r}\right);&& \label{eq:y:nestt}\\
		\lambda^{r+1}_{j}&=\lambda_{j}^{r}, \quad x^{r+1}_{j} = z^r, \quad \forall~j\ne i_{r}; \label{eq:x_i:nestt2}&&\\
		\mbox{Update $z$:}  \quad \quad  z^{r+1}&=\arg\min_{z\in Z}L(\{x_{ i}^{r+1}\}, z;  \lambda^{r}).\label{eq:z:nestt}
		\vspace{-0.4cm}
		\end{flalign}
		\ENDFOR
		\STATE {\bfseries Output:} $(z^m, x^m, \lambda^m) $ where $m$ randomly picked from $\{1,2,\cdots, R\}$.
	\end{algorithmic}
\end{algorithm}
We remark that NESTT-G is related to the popular ADMM method for {\it convex} optimization \cite{BoydADMMsurvey2011}.  However our particular update schedule (randomly picking $(x_i, \lambda_i)$ plus deterministic updating $z$), combined with the special $x$-step (minimizing an approximation of $L(\cdot)$ evaluated at {a different} block variable $z$) is not known before. These features are critical in our following rate analysis.

\subsection{Convergence Analysis.} To proceed, let us define $r(j)$ as the last iteration in which  the $j$th block is picked before iteration $r+1$. i.e.$r(j) := \max\{t\mid t<r+1, j = i(t)\}.$
Define $y_j^{r}: =z^{r(j)}$ if $j\ne i_{r}$, and $y_{i_{r}}^r=z^r$. Define the filtration $\cF^{r}$ as the $\sigma$-field generated by $\{i(t)\}_{t=1}^{r-1}$.

A few important observations are in order. Combining the $(x,z)$ updates \eqref{eq:x_i:nestt1} -- \eqref{eq:z:nestt}, we have {\small 
	\begin{subequations}
		\begin{align}
		&x^{r+1}_{q} = z^r -\frac{1}{\alpha_{q}\eta_{q}} (\lambda^r_{q} +\frac{1}{N}\nabla g_{q}(z^r)), ~\frac{1}{N}\nabla g_{q}(z^r) + \lambda^r_{q} +\alpha_{q}\eta_{q}(x^{r+1}_{q}-z^r) = 0, \; {\mbox \rm with}\; q = i_{r}\label{key:x:ir}\\
		&\lambda^{r+1}_{i_{r}} = -\frac{1}{N}\nabla g_{i_{r}}(z^r), \; \lambda^{r+1}_j = -\frac{1}{N}\nabla g_j(z^{r(j)}), \; \forall~j\ne i_{r}, \; \Rightarrow \lambda^{r+1}_{i} = -\frac{1}{N}\nabla g_i(y^r_i), \; \forall~i \label{key:lambda}\\
		&x^{r+1}_{j} \stackrel{\eqref{eq:x_i:nestt2}}=  z^r \stackrel{\eqref{key:lambda}}= z^r -\frac{1}{\alpha_{j}\eta_{j}} (\lambda^r_{j} +\frac{1}{N}\nabla g_{j}(z^{r(j)})) , \; \forall~j\ne i_{r}\label{key:x:j}.
		\end{align}
\end{subequations}}
\!\! The key here  is that the dual variables serve as the ``memory" for the past gradients of $g_i$'s. To proceed, we first construct a {\it potential function} using an {\it upper bound} of $L(x,y;\lambda)$. Note that {\small
	\begin{align}
	&\frac{1}{N}g_j(x_j^{r+1}) +\langle \lambda^r_j, x_j^{r+1} - z^r\rangle+\frac{\eta_j}{2}\|x^{r+1}_j-z^{r}\|^2 =  \frac{1}{N}g_j(z^r), \; \forall~j \ne i_{r} \\
	&\frac{1}{N}g_{i_{r}}(x_{i_{r}}^{r+1}) +\langle \lambda^r_{i_{r}}, x_{i_{r}}^{r+1} - z^r\rangle+\frac{\eta_i}{2}\|x^{r+1}_{i_{r}}-z^{r}\|^2\nonumber\\
	& \stackrel{\rm (i)}\le \frac{1}{N}g_{i_{r}}(z^r) +\frac{\eta_{i_{r}}+L_{i_{r}}/N}{2}\|x^{r+1}_{i_{r}}-z^{r}\|^2 \nonumber\\
	&\stackrel{\rm (ii)}= \frac{1}{N}g_{i_{r}}(z^r) +\frac{\eta_{i_{r}}+L_{i_{r}}/N}{2 (\alpha_{i_{r}}\eta_{i_{r}})^2}\|1/N (\nabla g_{i_{r}}(y_{i_{r}}^{r-1})-\nabla g_{i_{r}}(z^r))\|^2 
	\end{align}}
\!\! where ${\rm (i)}$ uses \eqref{key:lambda} and applies the descent lemma on the function $1/N g_i(\cdot)$; in ${\rm (ii)}$ we have used \eqref{eq:y:nestt} and \eqref{key:lambda}.
Since each $i$ is picked with probability $p_i$, we have
\begin{align*}
&\mathbb{E}_{i_{r}} [L(x^{r+1}, z^r; \lambda^r)\mid \cF^r] \nonumber\\
&\le \sum_{i=1}^{N}\frac{1}{N} g_i(z^r) +\sum_{i=1}^{N}\frac{p_i (\eta_{i}+L_{i}/N)}{2 (\alpha_{i}\eta_{i})^2}\|1/N(\nabla g_{i}(y_i^{r-1})-\nabla g_{i}(z^r))\|^2 + g_0(z^r) + h(z^r)\nonumber\\
&\le \sum_{i=1}^{N}\frac{1}{N} g_i(z^r) +\sum_{i=1}^{N}\frac{3p_i \eta_{i}}{(\alpha_{i}\eta_{i})^2}\|1/N(\nabla g_{i}(y_i^{r-1})-\nabla g_{i}(z^r))\|^2  + g_0(z^r) + h(z^r) : = Q^r,
\end{align*}
where in the last inequality we have used Assumption [A-(c)]. In the following, we will use $\mathbb{E}_{\cF^{r}}[Q^r]$ as the potential function, and show that it decreases at each iteration. 
\begin{lemma}\label{lem:nesttg:des}
	Suppose Assumption A holds, and pick 
	\begin{align}\label{eq:p:eta}
	\alpha_i= p_i =\beta \eta_i, \; \mbox{where}\; \beta:= \frac{1}{\sum_{i=1}^N\eta_i}, \quad \mbox{and} \quad {\eta_i\geq \frac{9L_i}{ N p_i}}, \quad  i=1,\cdots N.
	\end{align}
	Then the following descent estimate holds true for NESTT-G{
		\begin{align}\label{eq:nesttg:descent}
		\mathbb{E}[Q^{r}-Q^{r-1} | \cF^{r-1}]&\le -\frac{\sum_{i=1}^N\eta_i}{8}\mathbb{E}_{z^r}\|z^{r}-z^{r-1}\|^2 - \sum_{i=1}^N\frac{1}{{2\eta_i}}\|\frac{1}{N}(\nabla g_i(z^{r-1})-\nabla g_i(y_i^{r-2}))\|^2.
		\end{align}}
\end{lemma}
\noindent{\bf Sublinear Convergence.} Define the optimality gap as the following:{\small
	\begin{align}\label{eq:gap:def}
	\mathbb{E}[G^r] := \mathbb{E}\left[\|\tilde{\nabla}_{1/\beta} f(z^r)\|^2\right]=\frac{1}{\beta^2}\mathbb{E}\left[\|z^r-\prox_h^{1/\beta}[z^r-\beta \nabla (g(z^r)+g_0(z^r)) ]\|^2\right].
	\end{align}}
\!\! Note that when $h, g_0\equiv0$, $\mathbb{E} [G^r]$ reduces to $E[\|\nabla g(z^r)\|^2]$. We have the following result.
\begin{theorem}\label{thm:rate:nestt-g}
	Suppose Assumption A holds, and pick  (for $i=1,\cdots, N$)
	\begin{align}\label{eq:parameter:nesttg}
	\alpha_i= p_i=\frac{\sqrt{L_i/N}}{\sum_{i=1}^{N}\sqrt{L_i/N}}, \;  \eta_i = { 3} \left(\sum_{i=1}^{N}\sqrt{L_i/N}\right)\sqrt{L_i/N}, \;  \beta = \frac{1}{{3} (\sum_{i=1}^{N}\sqrt{L_i/N})^2}.
	\end{align}
	Then every limit point generated by NESTT-G is a stationary solution of problem \eqref{eq:distributed}. 
	Further, {\small
		\begin{align*}
		&1)~ \mathbb{E}[G^m]\le {\frac{80}{3}}\Big(\sum_{i=1}^N\sqrt{L_i/N}\Big)^2\frac{\mathbb{E}[{Q^1-Q^{R+1}} ]}{R};\\
		&2)~\mathbb{E}[G^m] + \mathbb{E}\left[\sum_{i=1}^N{3 \eta^2_i }\left\|x^{m}_i-z^{m-1}\right\|^2\right]\leq \frac{80}{3} \left(\sum_{i=1}^N\sqrt{L_i/N}\right)^2\frac{\mathbb{E}[{Q^1-Q^{R+1}} ]}{R}.
		\end{align*}}
\end{theorem}
{Note that Part (1) is useful in the {\it centralized} finite-sum minimization setting, as it shows the sublinear convergence of NESTT-G, measured only by the primal optimality gap evaluated at $z^r$. Meanwhile, part (2) is useful in the {\it distributed} setting, as it also shows that the expected constraint violation, which measures the consensus among agents,  shrinks in the same order. } We also comment that the above result suggests that to achieve an $\epsilon$-stationary solution, the NESTT-G requires about $\mathcal{O}\left(\bigg(\sum_{i=1}^{N}\sqrt{L_i/N}\bigg)^2/\epsilon\right)$ number of gradient evaluations (for simplicity we have ignored an additive $N$ factor for evaluating the gradient of the entire function at the initial step of the algorithm). 

It is interesting to observe that our choice of $p_i$ is proportional to the {\it square root} of the Lipschitz constant of each component function, rather than to $L_i$. Because of such choice of the sampling probability,  the derived convergence rate has a mild dependency on $N$ and $L_i$'s. Compared with the conventional gradient-based methods, our scaling can be up to $N$ times better. Detailed discussion and comparison will be given in Section \ref{sec:discuss}.

Note that similar sublinear convergence rates can be obtained for the case $\alpha_i =1$ for all $i$ (with different scaling constants). However due to space limitation, we will not present those results here. 

\subsection{Linear Convergence.}
In this section we show that the NESTT-G is capable of linear convergence for a family of nonconvex quadratic problems, which has important applications, for example in high-dimensional statistical learning \cite{loh12}. To proceed, we will assume the following. 
\begin{enumerate}
	\item [B-(a)] Each function $g_i(z)$ is a quadratic function of the form $g_i(z) = 1/2 z^T A_i z+ \langle b, z\rangle$, where $A_i$ is a symmetric matrix but not necessarily positive semidefinite;
	\item [B-(b)] The feasible set $Z$ is a closed compact polyhedral set;  
	\item [B-(c)] The nonsmooth function $p(z)=\mu\|z\|_1$, for some $\mu\ge 0$.
\end{enumerate}
Our linear convergence result is based upon certain error bound condition around the stationary solutions set, which has been shown in \cite{Luo92linear_convergence} for smooth quadratic problems and has been extended to including $\ell_1$ penalty in \cite[Theorem 4]{tseng09coordiate}.

\begin{lemma}\label{lm:eb}
	Suppose Assumptions A and B hold. Let  $Z^*$ denotes the set of stationary solutions of problem \eqref{eq:original}, and $\dist(z,Z^*):=\min_{u\in Z^*}\|z-u\|$. Then we have the following 
	\begin{enumerate}
		\item ({\bf Error Bound Condition}) For any $\xi\ge \min_z f(z)$,  exists a positive scalar $\tau$ such that the following error bound holds
		\begin{equation}\label{eq:primaleb}
		\dist(z,Z^*)\le \tau \|\tilde\nabla_{1/\beta} f(z)\|
		\end{equation}
		for all $z\in (Z\cap {\rm dom} ~h)$ and $z\in \{z: f(z)\le \xi \}$.
		\item ({\bf Separation of Isocost Surfaces}) There exists a scalar $\delta>0$ such that
		\begin{align}
		\|z - v\| \ge \delta\quad \mbox{whenever}\quad z\in Z^*, v\in Z^*, f(z) \ne f(v). 
		\end{align}
	\end{enumerate}
\end{lemma}
We note that the first statement holds true largely due to \cite[Theorem 4]{tseng09coordiate}, and the second statement holds true due to \cite[Lemma 2.1]{luo92}; see detailed discussion after \cite[Assumption 2]{tseng09coordiate}. Here the only difference with the statement  \cite[Theorem 4]{tseng09coordiate}  is that the error bound condition \eqref{eq:primaleb} holds true {\it globally}. This is by the assumption that $Z$ is a compact set. The proof will be provided in the Appendix.

Utilizing the above result, we have the following linear convergence claim.  
\begin{theorem}\label{thm:lin:conv:nesttg}
	Suppose that Assumptions A, B are satisfied. Then the sequence $\{\mathbb{E}[Q^{r+1}]\}_{r=1}^{\infty}$ converges $Q$-linearly \footnote{{A sequence
			$\{x^r\}$ is said to converge $Q$-linearly to some  $\bar{x}$ if
			$\lim\sup_r\|x^{r+1}-\bar{x}\|/\|x^{r}-\bar{x}\|\le \rho$, where
			$\rho\in(0,1)$ is some constant; cf \cite{tseng09coordiate} and references therein. }} to some $Q^* = f(z^*)$, where $z^*$ is a stationary solution for problem \eqref{eq:original}. That is, there exists a finite $\bar{r}>0$, $\rho\in (0,1)$ such that for all $r\ge \bar{r}$, $\mathbb{E}[Q^{r+1}-Q^*]{\le \rho} \mathbb{E}[Q^{r}-Q^*]$.
\end{theorem}

Linear convergence of this type for problems satisfying Assumption B has  been shown for (deterministic) proximal gradient based methods \cite[Theorem 2, 3] {tseng09coordiate}. To the best of our knowledge, this is the first result that shows the same linear convergence for a stochastic and distributed algorithm. There has been some recent works showing linear convergence for nonconvex problems satisfying certain quadratic growth condition {\cite{Karimi14,reddi16fast,zhu16VR}}. However the problems considered in {\cite{Karimi14,reddi16fast,zhu16VR}} are smooth unconstrained problems, and every stationary point is a global minimum, therefore they do no cover our nonconvex quadratic problems, whose stationary solutions are not global minimizers. 

\section{The NESTT-E Algorithm}
\subsection{ Algorithm Description}
In this section, we present a variant of NESTT-G, which is named NESTT with Exact minimization (NESTT-E).  
Our motivation is the following. First, in NESTT-G every agent should update its local variable at every iteration [cf. \eqref{eq:x_i:nestt1} or \eqref{eq:x_i:nestt2}]. In practice this may not be possible, for example at any given time a few agents can be in the {\it sleeping mode} so they cannot  perform \eqref{eq:x_i:nestt2}. Second, in the distributed setting it has been generally observed (e.g., see \cite[Section V]{chang14distributed}) that performing exact minimization (whenever possible) instead of taking the gradient steps for local problems can significantly speed up the algorithm. The NESTT-E algorithm to be presented in this section is designed to address these issues. 
\begin{algorithm}[tb]
	\caption{NESTT-E Algorithm}
	\label{alg:nestte}
	\begin{algorithmic}[1]
		\FOR{$r=1$ {\bfseries to} $R$}
		\STATE    Update $z$ by minimizing the augmented Lagrangian:
		\begin{flalign}
		z^{r+1}={\rm arg}\min_{z}\; L(x^{r}, z;  \lambda^{r}).&&\label{eq:z_update:exact}
		\end{flalign}
		\vspace{-.5cm}
		\STATE Randomly pick $i_{r}\in\{1,2,\cdots N\}$ with probability $p_{i_{r}}$:
		\begin{flalign}
		x^{r+1}_{i_{r}}&=\arg\!\min_{x_{i_{r}}} U_{i_{r}}(x_{i_{r}},z^{r+1};\lambda^r_{i_{r}});&& \label{eq:xi_update:exact}\\
		\lambda^{r+1}_{i_{r}}&=\lambda_{i_{r}}^{r}+\alpha_{i_{r}}\eta_{i_{r}}\left(x^{r+1}_{i_{r}}-z^{r+1}\right);&& \label{eq:y_update:exact}\\
		x^{r+1}_j &= x^{r}_j,\quad \lambda^{r+1}_j = \lambda^{r}_j \quad \forall~j\neq i_{r}.&&\label{eq:xi_update:rest}
		\end{flalign}
		\ENDFOR
		\STATE {\bfseries Output:} $(z^m, x^m, \lambda^m) $ where $m$ randomly picked from $\{1,2,\cdots, R\}$.
	\end{algorithmic}
\end{algorithm}
To proceed, let us define a new function as follows:
{\small
	\begin{align*}
	U(x, z; \lambda) := \sum_{i=1}^{N} U_i(x_i,z;\lambda_i):=\sum_{i=1}^{N}\left(\frac{1}{N}g_i(x_i)+\langle \lambda_i, x_i-z\rangle+\frac{\alpha_i\eta_i}{2}\|x_i-z\|^2\right).
	\end{align*}}
Note that if $\alpha_i=1$ for all $i$, then the $L(x, z; \lambda)  = U(x, z; \lambda) +p(z)+h(z) $. The algorithm details are presented in Algorithm 2.
The algorithm proceeds as follows.  At each iteration the cluster center minimizes $L(x, z; \lambda)$ with respect to $z$. Then the updated $z$ is sent to a randomly selected agent $i_{r}\in\{1,2,\cdots N\}$, who  minimizes $U(x, z; \lambda)$ w.r.t. its local variable $x_{i_{r}}$, followed by a dual ascent step for $\lambda_{i_{r}}$. 

\subsection{Convergence Analysis} 
We begin analyzing NESTT-E. The proof technique is quite different from that for NESTT-G, and it is based upon using the expected value of the {\it Augmented Lagrangian} function as the potential function.  For the ease of description we define the following quantities:
{
	\begin{align*}
	&w := (x, z, \lambda), \quad \beta:=\frac{1}{\sum_{i=1}^N\eta_i}, \quad c_i:= \frac{L^2_i}{\alpha_i\eta_i N^2}-\frac{\gamma_i}{2}+\frac{1-\alpha_i}{\alpha_i}\frac{L_i}{N}, \quad \alpha := \{\alpha_i\}_{i=1}^{N}.
	\end{align*}}
To measure the optimality of NESTT-E, define the {\it prox-gradient} of $L(x, z; \lambda)$ as:
\begin{align}
\tilde{\nabla}L(w)=\bigg[(z-\prox_{h}[z-\nabla_{z}(L(w)-h(z))]); \nabla_{x_1}L(w); \cdots;\nabla_{x_N} L(w)\bigg]\in \mathbb{R}^{(N+1) d}.
\end{align}
We define the optimality gap by adding to $\|\tilde{\nabla}L(w)\|^2$  the size of the constraint violation \cite{hong14nonconvex_admm}: 
{$$H(w^r):=\|\tilde{\nabla}L(w^r)\|^2+\sum_{i=1}^{N}\frac{L_i^2}{N^2}\|x^{r}_i-z^r\|^2.$$}
It can be verified that $H(w^r)\to 0$ implies that $w^r$ reaches a stationary solution for problem \eqref{eq:distributed}. 
{ We have the following theorem regarding the convergence properties of NESTT-E.}
\begin{theorem}\label{thm:sublin}
	Suppose Assumption A holds, and that $(\eta_i,\alpha_i)$ are chosen such that $c_i<0$ .  Then for some constant $\underline{f}$, we have
	$${\mathbb{E}[L(w^{r})]\ge \mathbb{E}[L(w^{r+1})]}\ge \underline{f}>-\infty, \quad \forall~r\ge 0.$$
	Further, almost surely every limit point of $\{w^r\}$ is a stationary solution of problem \eqref{eq:distributed}.
	Finally, for some function of $\alpha$ denoted as $C(\alpha)={\sigma_1(\alpha)}/{\sigma_2(\alpha)}$, we have the following:
	{
		\begin{align}
		\mathbb{E}[H(w^m)]\leq \frac{C(\alpha) \mathbb{E}[L(w^1)-L(w^{R+1})]}{R},
		\end{align}}
	where $\sigma_1:=\max(\hat{\sigma}_1(\alpha),\tilde{\sigma}_1)$ and $\sigma_2:=\max(\hat{\sigma}_2(\alpha),\tilde{\sigma}_2)$, and these constants are given by{
		\begin{align*}
		\hat{\sigma}_1(\alpha) &= \max_i\left\{4\left(\frac{L_i^2}{N^2}+\eta_i^2+\left(\frac{1}{\alpha_i}-1\right)^2\frac{L^2_i}{N^2}\right)+3\left(\frac{L_i^4}{\alpha_i \eta_i^2 N^4}+\frac{L_i^2}{N^2}\right)\right\},\nonumber\\
		\tilde{\sigma}_1&=\sum_{i=1}^N4\eta_i^2 + ({2}+\sum_{i=1}^N\eta_i+ L_0 )^2+3\sum_{i=1}^N\frac{L_i^2}{N^2},\nonumber\\
		\hat{\sigma}_2(\alpha) &= \max_i \left\{p_i\left(\frac{\gamma_i}{2}-\frac{L^2_i}{N^2\alpha_i\eta_i}-\frac{1-\alpha_i}{\alpha_i}\frac{L_i}{N}\right)\right\}, \quad \tilde{\sigma}_2 = \frac{\sum_{i=1}^N\eta_i-L_0}{2}.
		\end{align*}}
\end{theorem}
We remark that the above result shows the sublinear convergence of NESTT-E  to the set of stationary solutions. Note that $\gamma_i = \eta_i - L_i/N$, to satisfy $c_i<0$, a simple derivation yields 
{$$\eta_i > \frac{L_i \left((2-\alpha_i)+\sqrt{{(\alpha_i-2)^2} + {8 \alpha_i}}\right)}{ 2N \alpha_i}.$$}
Further, the above result characterizes the dependency of the rates on various parameters of the algorithm. For example, to see the effect of $\alpha$ on the convergence rate, let us set $p_i=\frac{L_i}{\sum_{i=1}^NL_i}$, and $\eta_i=3L_i/N$, and assume $L_0=0$, then consider two different choices of $\alpha$: $\widehat{\alpha}_i = 1, \; \forall~i$ and  $\widetilde{\alpha}_i = 4, \; \forall~i$. One can easily check that applying these different choices leads to following results:
{
	\begin{align*}
	C(\widehat{\alpha}) &= {49}\sum_{i=1}^NL_i/N, \quad \quad C(\widetilde{\alpha}) = 28\sum_{i=1}^NL_i/N.
	\end{align*}}
The key observation is that increasing $\alpha_i$'s reduces the constant in front of the rate. Hence, we expect that in practice larger $\alpha_i$'s will yield faster convergence. This phenomenon will be later confirmed by the numerical results.

Next let us briefly present the linear convergence of NESTT-E algorithm under Assumption B. The proof again utilizes the error bound condition in Lemma \ref{lm:eb}.
\begin{theorem}\label{thm:lin:conv:nestte}
	Suppose that Assumptions A, B are satisfied. Then the sequence $\{\mathbb{E}[L^{r+1}]\}_{r=1}^{\infty}$ converges $Q$-linearly  to some $L^* = f(z^*)$, where $z^*$ is a stationary solution for problem \eqref{eq:original}. That is, there exists a finite $\bar{r}>0$, $\rho\in (0,1)$ such that for all $r\ge \bar{r}$, $\mathbb{E}[L^{r+1}-L^*]\le\rho \mathbb{E}[L^{r}-L^*]$.
\end{theorem}

\section{Connections and Comparisons with Existing Works}\label{sec:discuss}

In this section we compare NESTT-G/E with a few existing algorithms in the literature. 
First, we present a somewhat surprising observation, that NESTT-G takes the same form as some well-known algorithms for {\it convex} finite-sum problems. 
To formally state such relation, we show in the following result that NESTT-G in fact admits a compact {\it primal-only} characterization. 
\begin{proposition}\label{prop:omega:close}
	The NESTT-G can be written into the following compact form:
	\begin{subequations}\label{eq:omega:close}{
			\begin{align}
			z^{r+1}&=\arg\min_z \; h(z)+ g_0(z)+\frac{1}{2\beta}\|z-u^{r+1}\|^2\\
			\mbox{\rm with}\quad u^{r+1}&:=z^r-\beta\Big(\frac{1}{{ N\alpha_{i_{r}}}}(\nabla g_{i_{r}}(z^r)-\nabla g_{i_{r}}(y_{i_{r}}^{r-1}))+\frac{1}{N}\sum_{i=1}^{N}\nabla g_i(y_i^{r-1})\Big).
			\end{align}}
	\end{subequations}
\end{proposition}
\! Based on this observation, the following comments are in order. 
\begin{enumerate}[leftmargin=*]
	\item[(1)] Suppose $h\equiv 0$, $g_0\equiv 0$ and $\alpha_{i} = 1$, $p_i=1/N$ for all $i$. Then \eqref{eq:omega:close}  takes the same form as the  SAG presented in \cite{Schmidt13}.  Further, when the component functions $g_i$'s are picked {\it cyclically } in a Gauss-Seidel manner, the iteration \eqref{eq:omega:close} takes the same form as the  IAG algorithm \cite{Blatt07}.
	\item[(2)] Suppose $h\ne 0$ and {$g_0\ne 0$}, and  $\alpha_{i} =p_i= 1/N$ for all $i$. Then \eqref{eq:omega:close}  is the same as the SAGA algorithm \cite{Defazio14}, which is design for optimizing convex nonsmooth finite sum problems.  \\
\end{enumerate}
Note that SAG/SAGA/IAG are all designed for convex problems. Through the lens of primal-dual splitting, our work shows that they can be generalized to nonconvex nonsmooth problems as well. 

Secondly, NESTT-E is related to the proximal version of the nonconvex ADMM \cite[Algorithm 2]{hong14nonconvex_admm}. However, the introduction of $\alpha_i$'s is new, which can significantly improve the practical performance but complicates the analysis. Further, there has been no counterpart of the sublinear and linear convergence rate analysis for the stochastic version of \cite[Algorithm 2]{hong14nonconvex_admm}. 

Thirdly, we note that a recent paper \cite{reddi16fast} has shown that SAGA works for smooth and unconstrained nonconvex problem. Suppose that $h\equiv 0$, {$g_0\ne 0$}, $L_i=L_j,\; \forall~i,j$ and  $\alpha_{i} = p_i = 1/N$, the authors show that SAGA achieves $\epsilon$-stationarity using $\mathcal{O}(N^{2/3}(\sum_{i=1}^{N}L_i/N)/\epsilon)$ gradient evaluations. Compared with GD, which achieves  $\epsilon$-stationarity using $\mathcal{O}(\sum_{i=1}^{N}L_i/\epsilon)$ gradient evaluations in the worse case (in the sense that $\sum_{i=1}^{N}L_i/N = L$), the rate in \cite{reddi16fast} is $\mathcal{O}(N^{1/3})$ times better. {However, the algorithm in \cite{reddi16fast} is different from NESTT-G in two aspects: 1) it does not generalize to the nonsmooth constrained problem \eqref{eq:original}; 2) it samples two component functions at each iteration, while NESTT-G only samples once.} Further, the analysis and the scaling are derived for the case of uniform $L_i$'s, therefore it is not clear how the algorithm and the rates can be adapted for the non-uniform case. On the other hand, our NESTT works for the general nonsmooth constrained setting. The non-uniform sampling used in NESTT-G is well-suited for problems with non-uniform $L_i$'s, and our scaling can be up to $N$ times better than GD (or its proximal version) in the worst case. {Note that problems with non-uniform $L_i$'s for the component functions are common in applications such as sparse optimization and signal processing.  For example in LASSO problem the data matrix is often normalized by feature (or ``column-normalized" \cite{negahban2012}), therefore the $\ell_2$ norm of each row of the data matrix (which corresponds to the Lipschitz constant for each component function) can be dramatically different. }

In Table \ref{tab:grad} we list the comparison of the number of gradient evaluations for NESTT-G and GD, in the worst case (in the sense that $\sum_{i=1}^{N}L_i/N = L$). For simplicity, we omitted an additive constant of $\mathcal{O}(N)$ for computing the initial gradients. 

{
	\begin{table}[]
		\centering
		\caption{Comparison of \# of gradient evaluations for NESTT-G and GD in the worst case}
		\label{tab:grad}
		\begin{tabular}{lcc}
			\cline{2-3}
			& {\bf NESTT-G }                              & {\bf GD}                            \\ \hline
			{\bf \# of Gradient Evaluations}                                                                                          & $\mathcal{O}\left((\sum_{i=1}^N\sqrt{L_i/N})^2/\epsilon\right)$ & $\mathcal{O}\left(\sum_{i=1}^N L_i/\epsilon\right)$                          \\ \hline
			{\bf Case I}:  $L_i=1, ~\forall i$                                                                                        & $\mathcal{O}(N/\epsilon)$                     & $\mathcal{O}(N/\epsilon)$                           \\ \hline
			\begin{tabular}[c]{@{}l@{}} {\bf Case II} : $\mathcal{O}(N^{2/3})$ terms with $L_i=N^{2/3}$\\ the rest with $L_i=1$\end{tabular} & $\mathcal{O}(N/\epsilon)$                     & \multicolumn{1}{l}{$\mathcal{O}(N^{4/3}/\epsilon)$} \\ \hline
			\begin{tabular}[c]{@{}l@{}}{\bf Case II} : $\mathcal{O}(\sqrt{N})$ terms with $L_i=N$\\ the rest with $L_i=1$\end{tabular} & $\mathcal{O}(N/\epsilon)$                     & \multicolumn{1}{l}{$\mathcal{O}(N^{3/2}/\epsilon)$} \\ \hline
			\begin{tabular}[c]{@{}l@{}} {\bf Case IV} : $\mathcal{O}(1)$ terms with $L_i=N^2$\\ the rest with $L_i=1$\end{tabular} & $\mathcal{O}(N/\epsilon)$                     & \multicolumn{1}{l}{$\mathcal{O}(N^{2}/\epsilon)$} \\ \hline
		\end{tabular}
\end{table}}

\section{Numerical Results}
In this section we evaluate the performance of NESTT. 
Consider the high dimensional regression problem with noisy observation \cite{loh12}, where $M$ observations are generated by $y = X \nu + \epsilon$. Here $y\in\mathbb{R}^M$ is the observed data sample; $X\in\mathbb{R}^{M\times P}$ is the covariate matrix;  $\nu\in\mathbb{R}^{P}$ is the ground truth, and $\epsilon\in\mathbb{R}^{M}$ is the noise. Suppose that the covariate matrix is not perfectly known, i.e., we observe $A=X+W$ where  $W\in \mathbb{R}^{M\times P}$ is the noise matrix with known covariance matrix $\Sigma_W$. Let us define $\hat{\Gamma}:=1/M(A^\top A)-\Sigma_W$, and $\hat{\gamma}:=1/M (A^\top y)$. To estimate the ground truth $\nu$, let us consider the following (nonconvex) optimization problem posed in  \cite[problem (2.4)]{loh12} (where $R>0$ controls sparsity):
\begin{align}\label{eq:l1}
\min_z \; z^\top\hat{\Gamma}z-\hat{\gamma}z \quad \st \quad\|z\|_1\le R.
\end{align}
Due to the existence of noise, $\hat{\Gamma}$ is not positive semidefinite hence the problem is not convex. Note that this problem satisfies Assumption A-- B, then by Theorem \ref{thm:lin:conv:nesttg} NESTT-G converges Q-linearly. 

To test the performance of the proposed algorithm, we generate the problem following similar setups as \cite{loh12}. Let $X=(X_1;\cdots, X_N)\in \mathbb{R}^{M\times P}$ with $\sum_{i} N_i = M$ and each $X_i\in\mathbb{R}^{N_i\times P}$ corresponds to $N_i$ data points, and it is generated from i.i.d Gaussian. Here $N_i$ represents the size of each mini-batch of samples. Generate the observations $y_i = X_i \times \nu^* + \epsilon_i\in\mathbb{R}^{N_i}$,  where $\nu^*$ is a $K$-sparse vector to be estimated, and $\epsilon_i\in\mathbb{R}^{N_i}$ is the random noise. Let $W =[W_1; \cdots; W_N] $, with $W_i\in\mathbb{R}^{N_i\times P}$ generated with i.i.d Gaussian. Therefore we have $z^\top\hat{\Gamma}z = \frac{1}{N}\sum_{i=1}^{N}\frac{N}{M}z^\top \left(X^\top_i X_i - W_i^{\top}W_i\right)z$. We set $M=100.000$, $P=5000$, $N = 50$, $K = {22}\approx\sqrt{P}$,and $R=\|\nu^*\|_1$. In simulation, we perform a mini-batch version of the algorithms, meaning that we split the data matrix $A$ and labels $y$ into $M$ submatrices and store them into different nodes. Therefore, in each node we have $\hat{\Gamma}_i\in\mathbb{R}^{n_i\times p}$, and $\hat{\gamma_i}\in\mathbb{R}^{p}$ such that $\sum_{i=1}^M n_i=n$. Here we set $M=30$. 
We implement NESTT-G/E, the SGD, and the nonconvex SAGA proposed in \cite{reddi16fast} with stepsize $\beta = \frac{1}{3 L_{\max} N^{2/3}}$ (with $L_{\max} := \max_{i} L_i$). Note that the SAGA proposed in \cite{reddi16fast} {\it only} works for the unconstrained problems with uniform $L_i$, therefore when applied to \eqref{eq:l1} it is {\it not} guaranteed to converge. Here we only include it for comparison purposes. 

\begin{figure}
	\centering
	\begin{subfigure}[b]{0.4\textwidth}
		\includegraphics[width=\textwidth]{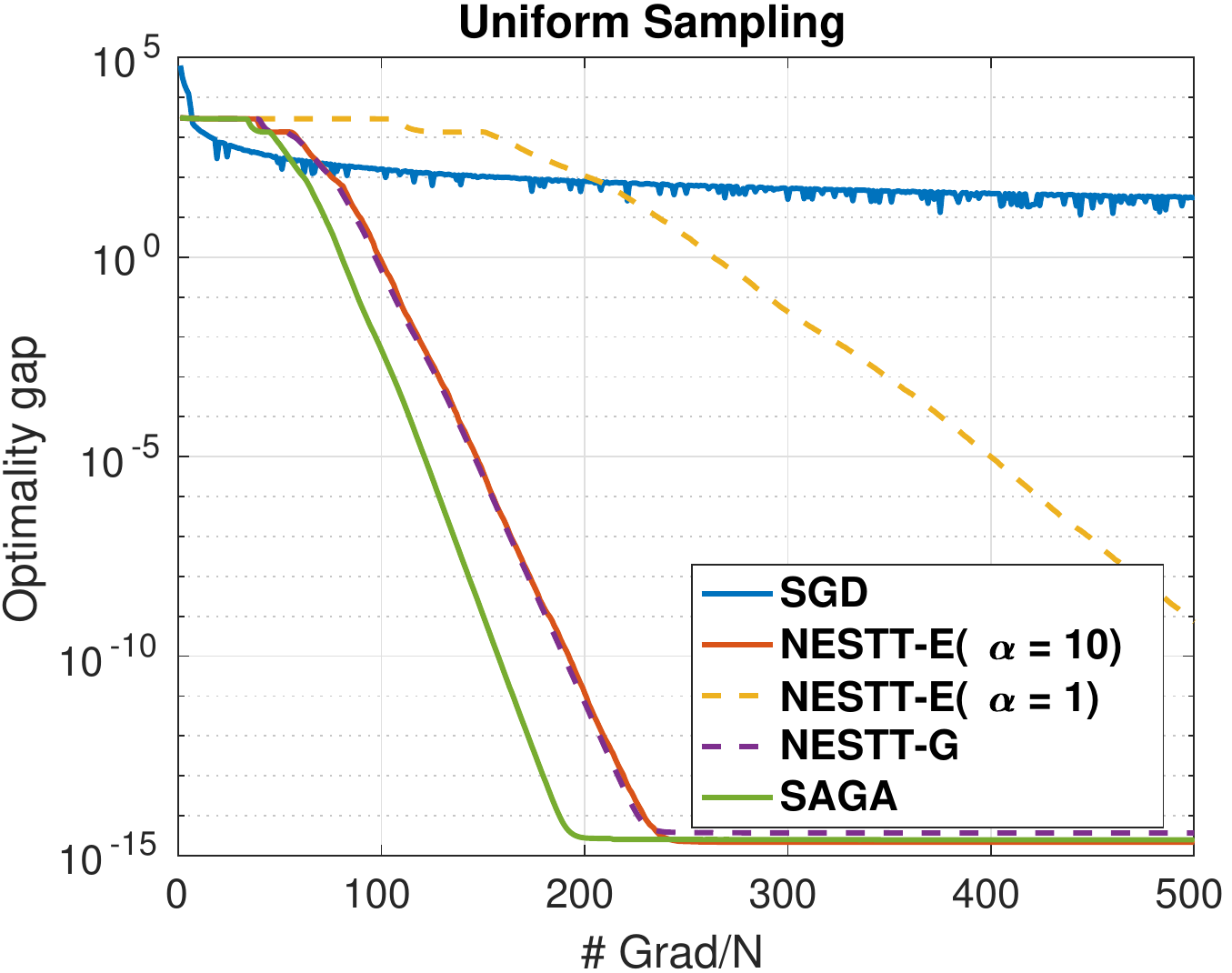}
		\vspace{-0.5cm}
		\label{fig:spca}
	\end{subfigure}
	\hfill
	\begin{subfigure}[b]{0.4\textwidth}
		\includegraphics[width=\textwidth]{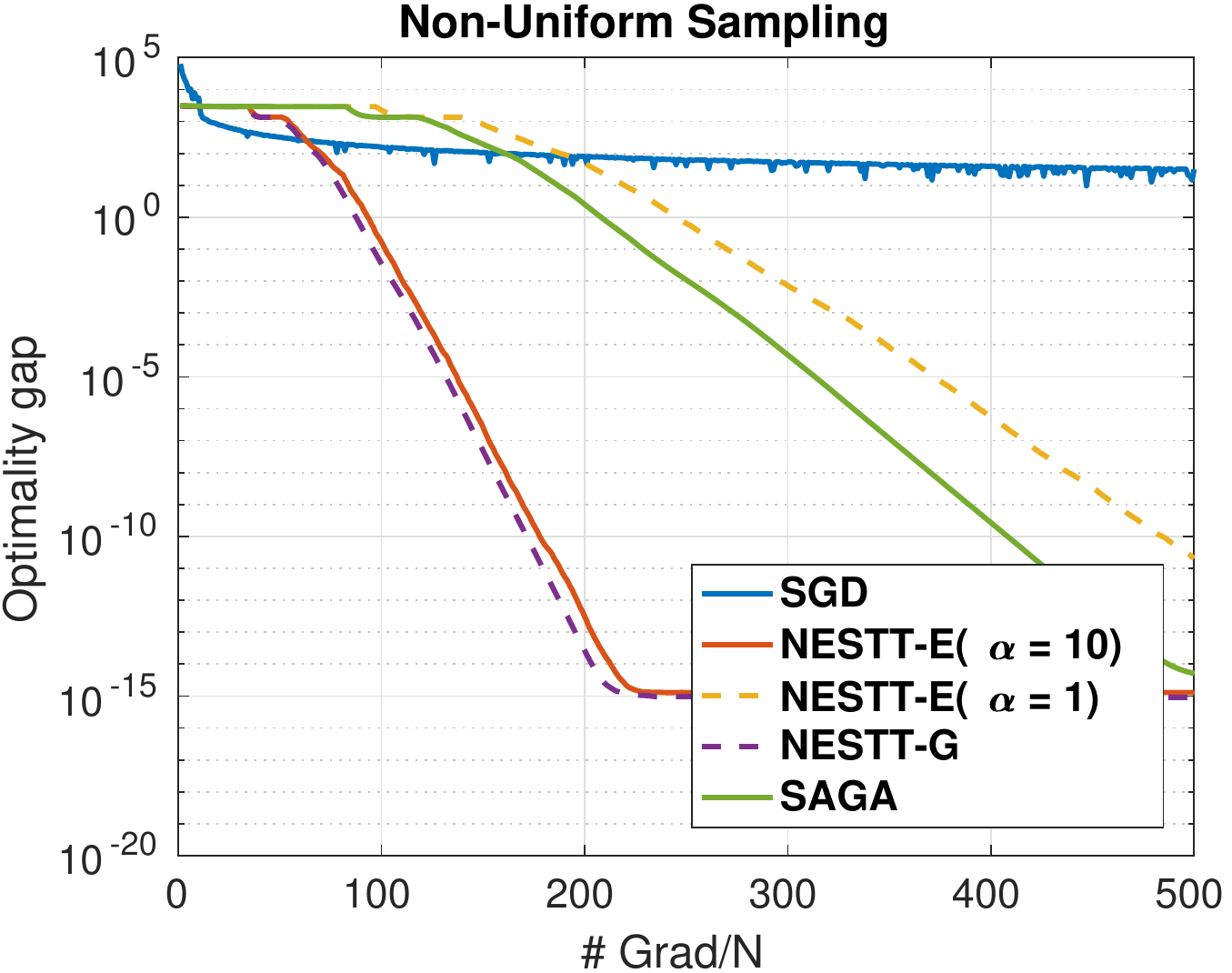}
		\vspace{-0.5cm}
		\label{fig:reg}
	\end{subfigure}
	\vspace{-0.1cm}
	\caption{{ \footnotesize Comparison of NESTT-G/E, SAGA, SGD on problem \eqref{eq:l1}. The $x$-axis denotes the number of passes of the dataset. Left: Uniform Sampling $p_i=1/N$; Right: Non-uniform Sampling ($p_i=\frac{\sqrt{L_i/N}}{\sum_{i=1}^N\sqrt{L_i/N}}$). }}
	\vspace{-0.2cm}
	\label{fig:1}      
\end{figure}

{In Fig. \ref{fig:1} we compare different algorithms in terms of the gap $\|\tilde{\nabla}_{1/\beta} f(z^r)\|^2$.} { In the left figure we consider the problem with $N_i=N_j$ for all $i,j$, and we show performance of the proposed algorithms with uniform sampling  (i.e., the probability of picking $i$th block is $p_i=1/N$). On the right one we consider problems in which approximately half of the component functions have twice the size of $L_i$'s as the rest, and consider the non-uniform sampling ($p_i=\sqrt{L_i/N}/{\sum_{i=1}^N\sqrt{L_i/N}}$). } Clearly in both cases the proposed algorithms perform quite well. Furthermore, it is clear that the NESTT-E performs well with large $\alpha:=\{\alpha_i\}_{i=1}^N$, which confirms our theoretical rate analysis. Also it is worth mentioning that when the $N_i$'s are non-uniform, the proposed algorithms [NESTT-G and NESTT-E (with $\alpha =10$)] significantly outperform SAGA and SGD. 
In Table \ref{table:2} we further compare different algorithms when changing the number of component functions (i.e., the number of mini-batches $N$) while the rest of the setup is as above. We run each algorithm with $100$ passes over the dataset. Similarly as before, our algorithms perform well, while SAGA seems to be sensitive to the uniformity of the size of the mini-batch [note that there is no convergence guarantee for SAGA applied to the nonconvex constrained problem \eqref{eq:l1}].

\begin{table}[]
	\centering\small
	\caption{Optimality gap $\|\tilde{\nabla}_{1/\beta} f(z^r)\|^2$ for different algorithms, with $100$ passes of the datasets. }\label{table:2}
	\label{tab:gap}
	\begin{tabular}{lllllllll}
		\cline{2-9}
		& \multicolumn{2}{c}{\textbf{SGD}}                          & \multicolumn{2}{c}{\textbf{NESTT-E ($\alpha =10$)}}                      & \multicolumn{2}{c}{\textbf{NESTT-G}}                      & \multicolumn{2}{c}{\textbf{SAGA}}                         \\ \hline
		N  & \multicolumn{1}{c}{Uniform} & \multicolumn{1}{c}{Non-Uni} & \multicolumn{1}{c}{Uniform} & \multicolumn{1}{c}{Non-Uni} & \multicolumn{1}{c}{Uniform} & \multicolumn{1}{c}{Non-Uni} & \multicolumn{1}{c}{Uniform} & \multicolumn{1}{c}{Non-Uni} \\ \hline
		10 & 3.4054                      & 0.2265                      & 2.6E-16                     & 6.16E-19                    & 2.3E-21                     & 6.1E-24                     & 2.7E-17                     & 2.8022                      \\
		20 & 0.6370                      & 6.9087                      & 2.4E-9                      & 5.9E-9                      & 1.2E-10                     & 2.9E-11                     & 7.7E-7                      & 11.3435                     \\
		30 & 0.2260                      & 0.1639                      & 3.2E-6                      & 2.7E-6                      & 4.5E-7                      & 1.4E-7                      & 2.5E-5                      & 0.1253                      \\
		40 & 0.0574                      & 0.3193                      & 5.8E-4                      & 8.1E-5                      & 1.8E-5                      & 3.1E-5                      & 4.1E-5                      & 0.7385                      \\
		50 & 0.0154                      & 0.0409                      & 8.3E.-4                      & 7.1E-4                      & 1.2E-4                      & 2.7E-4                      & 2.5E-4                      & 3.3187                      \\ \hline
	\end{tabular}
\end{table}

\newpage
{\centering\large{\bf Appendix}}
\noindent\subsection{Some Key Properties of NESTT-G}\label{sec:key:nesttg}
To facilitate the following derivation, in this section we collect some key properties of NESTT-G.

First, from the optimality condition of the $x$ update we have 
\begin{subequations}
	\begin{align}
	&x^{r+1}_{i_{r}} = z^r -\frac{1}{\alpha_{i_{r}}\eta_{i_{r}}} \left(\lambda^r_{i_{r}} +\frac{1}{N}\nabla g_{i_{r}}(z^r)\right), \label{key:x:ir:1}\\
	&x^{r+1}_{j} \stackrel{\eqref{eq:x_i:nestt2}}=  z^r \stackrel{\eqref{key:lambda}}= z^r -\frac{1}{\alpha_{j}\eta_{j}} (\lambda^r_{j} +\frac{1}{N}\nabla g_{j}(z^{r(j)})) , \; \forall~j\ne i_{r}\label{key:x:j:1}. 
	\end{align}
\end{subequations}
Then using the update scheme of the $\lambda$ we can further obtain
\begin{subequations}
	\begin{align}
	&\lambda^{r+1}_{i_{r}} = -\frac{1}{N}\nabla g_{i_{r}}(z^r), \\
	&\lambda^{r+1}_j = -\frac{1}{N}\nabla g_j(z^{r(j)}), \; \forall~j\ne i_{r}.\label{key:lambda:1}
	\end{align}
\end{subequations}
Therefore, using the definition of $y^r_i$ we have the following compact forms
\begin{align}
\lambda^{r+1}_{i} & = -\frac{1}{N}\nabla g_i(y^r_i), \; i=1,\cdots, N\label{eq:lambda:compact}. \\
x^{r+1}_{i} & = z^r -\frac{1}{\alpha_{i}\eta_{i}} \left(\lambda^r_{i} +\frac{1}{N}\nabla g_{i}(y^r_i)\right) , \; i=1,\cdots, N\label{eq:x:compact}.
\end{align}

Second, let us look at the optimality condition for the $z$ update. The $z$-update \eqref{eq:z:nestt} is given by
\begin{align}
z^{r+1}&=\arg\min_{z}\; L(\{x_{ i}^{r+1}\}, z;  \lambda^{r})\nonumber\\
&=\arg\min_{z} \; \sum_{i=1}^{N}\left( \langle \lambda^r_i, x^{r+1}_i-z\rangle+\frac{\eta_i}{2}\|x^{r+1}_i-z\|^2\right) + g_0(z) + h(z).
\end{align}
Note that this problem is strongly convex because we have assumed that $\sum_{i=1}{\eta_i} >  3 L_0$; cf. Assumption [A-(c)].

Let us define 
\begin{align}
u^{r+1} &:= \frac{\sum_{i=1}^{N} \eta_i x^{r+1}_i +\sum_{i=1}^{N}{\lambda^{r}_i}}{\sum_{i=1}^{N}\eta_i}\nonumber\\
& = \frac{\sum_{i=1}^{N} \eta_i z^r-  \eta_{i_{r}} (z^r - x_{i_{r}}^{r+1})}{\sum_{i=1}^{N}\eta_i}  + \frac{\sum_{i=1}^{N}{\lambda^{r}_i}}{\sum_{i=1}^{N}\eta_i} \nonumber\\
&\stackrel{\eqref{key:x:ir:1}} = \frac{\sum_{i=1}^{N} \eta_i z^r-  \frac{\eta_{i_{r}}}{\alpha_{i_{r}}\eta_{i_{r}}} (\lambda^r_{i_{r}}+1/N\nabla g_{i_{r}} (z^r))}{\sum_{i=1}^{N}\eta_i}  + \frac{\sum_{i=1}^{N}{\lambda^{r}_i}}{\sum_{i=1}^{N}\eta_i} \nonumber\\
& 	\stackrel{\eqref{eq:lambda:compact}} = z^r - \frac{\frac{1}{\alpha_{i_{r}}} (-\nabla g_{i_{r}}(y^{r-1}_{i_{r}})+\nabla g_{i_{r}} (z^r))}{N\sum_{i=1}^{N}\eta_i}  - \frac{\sum_{i=1}^{N}{\nabla g_i (y^{r-1}_i)}}{N\sum_{i=1}^{N}\eta_i} \nonumber\\
&\stackrel{\rm (i)}= z^r - {\frac{\beta}{N\alpha_{i_{r}}} (-\nabla g_{i_{r}}(y^{r-1}_{i_{r}})+\nabla g_{i_{r}} (z^r))} - \frac{\beta\sum_{i=1}^{N}{\nabla g_i (y^{r-1}_i)}}{N} \label{eq:comp}\\
&\stackrel{\rm (ii)}:= z^r - \beta v^{r+1}_{i_{r}}\label{eq:u}
\end{align}
where in ${\rm (i)}$ we have defined $\beta:= 1/\sum_{i=1}^{N}\eta_i$; in ${\rm (ii)}$ we have defined 
\begin{align}
v_{i_{r}}^{r+1} := \frac{1}{N}\sum_{i=1}^N\nabla g_i(y_i^{r-1}) + \frac{1}{ \alpha_{i_{r}}}\left(-\frac{1}{N}\nabla g_{i_{r}}(y^{r-1}_{i_{r}})+\frac{1}{N}\nabla g_{i_{r}}(z^{r})\right)\label{eq:v_ir}.
\end{align}
Clearly if we pick $\alpha_i = p_i$ for all $i$, then we have
\begin{align}\label{eq:unbiased}
\mathbb{E}_{i_{r}}[u^{r+1} \mid \cF^{r}] = z^r -\frac{\beta}{N} \sum_{i=1}^{N}\nabla g_i(z^r).
\end{align}
Using the definition of $u^{r+1}$, it is easy to check that 
\begin{align}\label{eq:z:prox}
z^{r+1} &= \arg\min_{z}\; \frac{1}{2\beta}\|z-u^{r+1}\|^2 + h(z) + g_0(z)\nonumber\\
& = \prox_h^{1/\beta}[ u^{r+1} - \beta \nabla g_0(z^{r+1}) ].
\end{align}

The optimality condition for the $z$ subproblem is given by:
\begin{align}
z^{r+1}-u^{r+1}+ \beta \nabla g_0(z^{r+1})+\beta \xi^{r+1} =0 \label{eq:opt:z:nestt}
\end{align}
where, $\xi^{r+1}\in \partial h(z^{r+1})$ is a subgradient of  $h(z^{r+1})$.  
Using the definition of $v_{i_{r}}$ in \eqref{eq:v_ir}, we obtain
\begin{align}
z^{r+1}=z^r-\beta(v^{r+1}_{i_{r}}+ \nabla g_0(z^{r+1}) +\xi^{r+1})\label{eq:z:ex}.
\end{align}

Third, if $\alpha_i=p_i$, then we have:
\begin{align}
&\mathbb{E}_{i_{r}}\left[\left\| -\frac{\lambda_{i_{r}}^{r}+1/N\nabla g_{i_{r}}(z^{r})}{\alpha_{i_{r}}}+\frac{1}{N}\sum_{i=1}^{N}\nabla g_i(z^{r})-\sum_{i=1}^N \frac{1}{N}\nabla g_i(y_i^{r-1})\right\|^2\right]\nonumber\\
&\stackrel{(a)}= \text{Var}\left[-\frac{\lambda_{i_{r}}^{r}+1/N\nabla g_{i_{r}}(z^{r})}{\alpha_{i_{r}}}\right]\nonumber\\
&\stackrel{(b)}\leq  \sum_{i=1}^N\frac{1}{\alpha_i}\left\|\frac{1}{N}\nabla g_i(z^{r})-\frac{1}{N}\nabla g_i(y_i^{r-1})\right\|^2,\label{eq:bd:var}
\end{align}
where $(a)$ is true because whenever $\alpha_i= p_i$ for all $i$, then $$\mathbb{E}_{i_{r}}\left[ -\frac{\lambda_{i_{r}}^{r}+1/N\nabla g_{i_{r}}(z^{r})}{\alpha_{i_{r}}}\right]=\frac{1}{N}\sum_{i=1}^{N}\nabla g_i(z^{r})-\sum_{i=1}^N \frac{1}{N}\nabla g_i(y_i^{r-1});$$
The inequality in $(b)$ is true because for a random variable $x$ we have $\text{Var}(x)\leq \mathbb{E}[x^2]$.

\noindent\subsection{Proof of Lemma \ref{lem:nesttg:des}}

{\bf Step 1).}  Using the definition of potential function $Q^r$, we have:
\begin{align}
&\mathbb{E}[Q^{r}- Q^{r-1}\mid \cF^{r-1}]\nonumber\\
&=\mathbb{E}\left[\sum_{i=1}^N\frac{1}{N}\left(g_i(z^r)-g_i(z^{r-1})\right) + g_0(z^r) -g_0(z^{r-1})+h(z^r)-h(z^{r-1})\mid {\cF^{r-1}}\right]\nonumber\\
&+\mathbb{E}\left[\sum_{i=1}^N\frac{3p_i}{\alpha_i^2\eta_i}\left\|\frac{1}{N}\nabla g_i(z^r)-\frac{1}{N}\nabla g_i(y_i^{r-1})\right\|^2 - \frac{3p_i}{\alpha_i^2\eta_i}\left\| \frac{1}{N}\nabla g_i(z^{r-1})-\frac{1}{N}\nabla g_i(y_i^{r-2})\right\|^2\mid  {\cF^{r-1}}\right]\label{eq:p:diff}.
\end{align}

{\bf Step 2).}  The first term in \eqref{eq:p:diff} can be bounded as follows (omitting the subscript $\cF^{r}$). 
\begin{align}
&\mathbb{E}\left[\sum_{i=1}^N \frac{1}{N}\left(g_i(z^r)-g_i(z^{r-1})\right)+ g_0(z^r)-g_0(z^{r-1})+h(z^r)-h(z^{r-1})\mid  {\cF^{r-1}}\right]\nonumber\\
&\stackrel{\rm (i)}\leq \mathbb{E}\bigg[\frac{1}{N}\sum_{i=1}^N\langle \nabla g_i(z^{r-1}),z^r-z^{r-1}\rangle+ \langle \nabla g_0(z^{r-1}),z^r-z^{r-1}\rangle \nonumber\\
&\quad \quad  +\langle \xi^r, z^r-z^{r-1}\rangle +\frac{\sum_{i=1}^N{L_i/N + L_0}}{2}\|z^r-z^{r-1}\|^2\mid  {\cF^{r-1}}\bigg]\nonumber\\
&\stackrel{\rm (ii)}=\mathbb{E}\left[\left\langle \frac{1}{N}\sum_{i=1}^N \nabla g_i(z^{r-1})+\xi^r + \nabla g_0(z^{r})+\frac{1}{\beta} (z^r-z^{r-1}),z^r-z^{r-1}\right\rangle\mid  {\cF^{r-1}}\right]\nonumber\\
&\quad \quad -\left(\frac{1}{\beta}-\frac{\sum_{i=1}^{N}L_i /N + 3 L_0}{2}\right)\mathbb{E}_{z^r}\|z^r-z^{r-1}\|^2\nonumber\\
&\stackrel{\eqref{eq:z:ex}}=\mathbb{E}\left[\left\langle\frac{1}{N}\sum_{i=1}^N \nabla g_i(z^{r-1})-v_{i(r-1)}^{r},z^r-z^{r-1}\right\rangle\mid  {\cF^{r-1}}\right]\nonumber\\
&\quad -\left(\frac{1}{\beta}-\frac{\sum_{i=1}^{N}L_i /N + 3 L_0}{2}\right)\mathbb{E}_{z^r}\|z^r-z^{r-1}\|^2\nonumber\\
&\stackrel{\rm (iii)}\leq \frac{1}{2\ell_1}\mathbb{E}\left[\left\|1/N\sum_{i=1}^N \nabla g_i(z^{r-1})-v_{i(r-1)}^{r}\right\|^2\mid \cF^{r-1}\right]+\frac{\ell_1}{2}\mathbb{E}_{z^r}\|z^r-z^{r-1}\|^2\nonumber\\
&\quad-\left(\frac{1}{\beta}-\frac{\sum_{i=1}^NL_i/N+3 L_0}{2}\right)\mathbb{E}_{z^r}\|z^r-z^{r-1}\|^2
\end{align}
where in ${\rm (i)}$ we have used the Lipschitz continuity of the gradients of $g_i$'s as well as the convexity of $h$; in ${\rm (ii)}$ we have used the fact that 
\begin{align}
\langle \nabla g_0(z^{r-1}),z^r-z^{r-1}\rangle  \le  \langle \nabla g_0(z^{r}),z^r-z^{r-1}\rangle + L_0\|z^{r}-z^{r-1}\|^2;
\end{align}
in ${\rm (iii)}$ we have applied the Young's inequality for some $\ell_1>0$.  

Choosing $\ell_1=\frac{1}{2\beta}$, we have:
\begin{align*}
&\frac{1}{2\ell_1}\mathbb{E}\left\|\frac{1}{N}\sum_{i=1}^N \nabla g_i(z^{r-1})-v_{i(r-1)}^{r}\right\|^2\nonumber\\
&\stackrel{\eqref{eq:v_ir}}=\beta\mathbb{E}\left[\left\| \frac{1}{N}\sum_{i=1}^{N}\nabla g_i(z^{r-1})-\frac{\lambda_{i(r-1)}^{r-1}+1/N\nabla g_{i(r-1)}(z^{r-1})}{\alpha_{i(r-1)}}-\sum_{i=1}^N \frac{1}{N}\nabla g_i(y_i^{r-2})\right\|^2\right]\\
&\stackrel{\eqref{eq:bd:var}}\leq \beta \sum_{i=1}^N\frac{1}{\alpha_i}\left\|\frac{1}{N}\nabla g_i(z^{r-1})-\frac{1}{N}\nabla g_i(y_i^{r-2})\right\|^2.
\end{align*}
Overall we have the following bound for the first term in \eqref{eq:p:diff}:
\begin{align}\label{eq:diff:p:1}
&\mathbb{E}\left[\sum_{i=1}^N\frac{1}{N}\left(g_i(z^r)-g_i(z^{r-1})\right)+ g_0(z^r) - g_0(z^{r-1})+h(z^r)-h(z^{r-1})\mid \cF^{r-1}\right]\\
&\leq \sum_{i=1}^N\frac{\beta }{\alpha_i}\left\|\frac{1}{N}\nabla g_i(z^{r-1})-\frac{1}{N}\nabla g_i(y_i^{r-2})\right\|^2-\left(\frac{3}{4\beta}-\frac{\sum_{i=1}^NL_i/N + 3L_0}{2}\right)\mathbb{E}_{z^r}\|z^r-z^{r-1}\|^2\nonumber.
\end{align} 

{
	{\bf Step 3).} We bound the second term in \eqref{eq:p:diff} in the following way:
	\begin{align}
	&\mathbb{E}\left[\|\nabla g_i(z^{r})-\nabla g_i(y_i^{r-1})\|^2 \mid \cF^{r-1}\right]\nonumber\\
	&=\mathbb{E}\left[\|\nabla g_i(z^{r})-\nabla g_i(y_i^{r-1})+\nabla g_i(z^{r-1})-\nabla g_i(z^{r-1})\|^2 \mid \cF^{r-1}\right]\nonumber\\
	&\stackrel{\rm (i)}\leq (1+\xi_i)\mathbb{E}_{z^r}\|\nabla g_i(z^{r})-\nabla g_i(z^{r-1})\|^2+\left(1+\frac{1}{\xi_i}\right)\mathbb{E}_{y_i^{r-1}}\|\nabla g_i(y_i^{r-1})-\nabla g_i(z^{r-1})\|^2\nonumber\\
	&\stackrel{\rm (ii)}= (1+\xi_i)\mathbb{E}_{z^r}\|\nabla g_i(z^{r})-\nabla g_i(z^{r-1})\|^2+(1-p_i)\left(1+\frac{1}{\xi_i}\right)\|\nabla g_i(y_i^{r-2})-\nabla g_i(z^{r-1})\|^2\label{eq:del:g_i}
	\end{align}
	where in ${\rm (i)}$ { we have used the fact that the randomness of $z^{r-1}$ comes from $i_{r-2}$, so fixing $\mathcal{F}^{r-1}$, $z^{r-1}$ is deterministic}; we have also applied the following inequality:
	\begin{align*}
	(a+b)^2\leq (1+\xi)a^2+(1+\frac{1}{\xi})b^2 \quad \forall ~ \xi>0.
	\end{align*}
	The equality ${\rm (ii)}$ is true because {the randomness of $y^{r-1}_i$ comes from $i_{r-1}$}, and for each  $i$ there is a probability $p_i$ such that $x^r_i$ is updated, so that  $\nabla g_i(y^{r-1}_i) = \nabla g_i(z^{r-1})$, otherwise $x_i$ is not updated so that $\nabla g_i(y^{r-1}_i) = \nabla g_i(y^{r-2}_i)$.
}

{\bf Step 4).} Applying  \eqref{eq:del:g_i} and set $\alpha_i=p_i$, the second part of \eqref{eq:p:diff} can be bounded as
\begin{align}
&\mathbb{E}\left[\sum_{i=1}^N\frac{3p_i}{\alpha_i^2\eta_i}\left\|\frac{1}{N}\nabla g_i(z^r)-\frac{1}{N}\nabla g_i(y_i^{r-1})\right\|^2 - \frac{3p_i}{\alpha_i^2\eta_i}\left\|\frac{1}{N}\nabla g_i(z^{r-1})-\frac{1}{N}\nabla g_i(y_i^{r-2})\right\|^2\mid \cF^{r-1}\right]\nonumber\\
&\le \sum_{i=1}^N\frac{3L_i^2}{\alpha_i\eta_i N^2}\left(1+\xi_i\right)\mathbb{E}_{z^r}\|z^{r}-z^{r-1}\|^2\nonumber\\
&+\frac{3}{\alpha_i\eta_i}\left((1-p_i)(1+\frac{1}{\xi_i})-1\right)\left\|\frac{1}{N}\nabla g_i(y_i^{r-2})-\frac{1}{N}\nabla g_i(z^{r-1})\right\|^2\label{eq:diff:p:2}.
\end{align}
Combining \eqref{eq:diff:p:1} and \eqref{eq:diff:p:2} eventually we have 
\begin{align}
&\mathbb{E}[Q^{r}- Q^{r-1}\mid \cF^r]\nonumber\\
&\le  \sum_{i=1}^N\left\{\frac{\beta}{\alpha_i} +\frac{3}{\alpha_i\eta_i}\left((1-p_i)(1+\frac{1}{\xi_i})-1\right)\right\}\left\|\frac{1}{N}\nabla g_i(z^{r-1})-\frac{1}{N}\nabla g_i(y_i^{r-2})\right\|^2\nonumber\\
& + \left\{-\frac{3}{4\beta}+\frac{\sum_{i=1}^NL_i/N + 3 L_0}{2}+\sum_{i=1}^N\frac{3L_i^2}{\alpha_i\eta_i N^2}\left(1+\xi_i\right)\right\}\mathbb{E}_{z^r}\|z^r-z^{r-1}\|^2.
\end{align}
Let us define $\{\tilde{c}_i\}$ and $\hat{c}$ as following:
\begin{align*}
\tilde{c}_i &= \frac{\beta}{\alpha_i} +\frac{3}{\alpha_i\eta_i}\left((1-p_i)(1+\frac{1}{\xi_i})-1\right)\\
\hat{c} &= -\frac{3}{4\beta}+\frac{\sum_{i=1}^NL_i/N + 3L_0}{2}+\sum_{i=1}^N\frac{3L_i^2}{\alpha_i\eta_i N^2}\left(1+\xi_i\right).
\end{align*}
In order to prove the lemma it is enough to show that  $\tilde{c}_i<-\frac{1}{2\eta_i} \; \forall~i$, and $\hat{c}<-\sum_{i=1}^N\frac{\eta_i}{8}$. Let us  pick 
\begin{align}\label{eq:parameter:select}
\alpha_i=p_i, \; \xi_i = {\frac{2}{p_i}}, \; p_i = \frac{\eta_i}{\sum_{i=1}^N\eta_i}.
\end{align}
Recall that {$\beta =\frac{1}{\sum_{i=1}^{N}\eta_i}$}.These values yield the following 
\begin{align*}
\tilde{c}_i = \frac{1}{\eta_i}-\frac{3}{\eta_i}\left(\frac{p_i+1}{2}\right)\leq \frac{1}{\eta_i}-\frac{3}{2\eta_i}=-\frac{1}{2\eta_i}< 0.
\end{align*}
To show that $\hat{c}\le -\sum_{i=1}^N\frac{\eta_i}{8}$ let us assume that $\eta_i=d_iL_i$ for some $d_i>0$. Note that by assumption we have 
$${\sum_{i=1}^{N}\eta_i\ge 3 L_0}.$$
{
	Therefore we have the following expression for $\hat{c}$:
	\begin{align*}
	\hat{c} &\le -\sum_{i=1}^N \frac{1}{4}d_iL_i+\frac{L_i}{2N } + \frac{3L_i}{p_id_i N^2}\left(1+\frac{2}{p_i}\right)\nonumber\\
	&<\sum_{i=1}^N \frac{L_i}{d_i}\left(-\frac{1}{4}d_i^2+\frac{d_i}{2 N}+\frac{9}{p_i^2 N^2}\right).
	\end{align*}
	As a result, to have $\hat{c}<-\sum_{i=1}^N\frac{\eta_i}{8}$, we need
	\begin{align}
	\frac{L_i}{d_i}\left(\frac{1}{4}d_i^2-\frac{d_i}{2 N}-\frac{9}{p_i^2 N^2}\right) \ge \frac{d_i L_i}{8}, \quad \forall~i.
	\end{align}
	Or equivalently
	\begin{align}
	\frac{1}{8}d_i^2-\frac{d_i}{2 N}-\frac{9}{p_i^2 N^2}\ge 0, \quad \forall~i.
	\end{align}
}
By finding the root of the above quadratic inequality, we need $d_i\geq \frac{9}{ N p_i}$,  which is equivalent to choosing the following parameters
\begin{align}\label{eq:eta}
\eta_i\geq\frac{9 L_i}{ N p_i}.
\end{align}

The lemma is proved. \QED
\subsection{\bf Proof of Theorem \ref{thm:rate:nestt-g}}
First, using the fact that $f(z)$ is lower bounded [cf. Assumption A-(a)], it is easy to verify that $\{Q^{r}\}$ is a bounded sequence. Denote its lower bound to be $\underbar{Q}$.  From Lemma \ref{lem:nesttg:des}, it is clear that $\{Q^r- \underbar{Q}\}$ is a nonnegative supermartingale. Apply the  Supermartigale Convergence Theorem \cite[Proposition 4.2]{bertsekas97} we conclude that $\{Q^r\}$ converges almost surely (a.s.), and that
\begin{align}\label{eq:delta:go:0}
\left\|\nabla g_i(z^{r-1})-\nabla g_i(y_i^{r-2})\right\|^2 \to 0, \quad \mathbb{E}_{z^r}\|z^{r}-z^{r-1}\| \to 0, \quad \mbox{a.s.}, \quad \forall~i. 
\end{align}
{The first inequality implies that $\|\lambda_{i_{r}}^r- \lambda_{i_{r}}^{r-1}\| \to 0$. Combining this with equation \eqref{eq:y:nestt} yields $\|x_{i_{r}}^r-z^{r-1}\|\to0$,  which further implies that $\|z^r- z^{r-1}\| \to 0$.} By utilizing \eqref{key:lambda} -- \eqref{key:x:j}, we can conclude that
\begin{align}
\quad \|x_i^{r}-x_i^{r-1}\| \to 0, \quad \|\lambda^r_i-\lambda^{r-1}_i\|\to 0, \quad \mbox{a.s.}, \quad \forall~i. 
\end{align}
That is, almost surely the successive differences of all the primal and dual variables go to zero. Then it is easy to show that every limit point of the sequence $(x^r, z^r, \lambda^r)$ converge to a stationary solution of problem \eqref{eq:distributed} (for example, see the argument in \cite[Theorem 2]{hong14nonconvex_admm}. Here we omit the full proof. 

{\bf Part 1).}  We bound the gap in the following way (where the expectation is taking over the nature history of the algorithm):
{
	\begin{align}\label{eq:Q:bound}
	&\mathbb{E}\left[\|z^r-\prox_h^{1/\beta}[z^r-\beta \nabla (g(z^r) + g_0(z^r))]\|^2\right]\nonumber\\
	&\stackrel{\rm(a)}= \mathbb{E}\left[\|z^r-z^{r+1}+\prox_h^{1/\beta}[u^{r+1} -\beta \nabla g_0(z^{r+1})]-\prox_h^{1/\beta}[z^r-\beta \nabla (g(z^r)+g_0(z^r)) ]\|^2\right]\nonumber\\
	&\stackrel{\rm(b)}\leq 3\mathbb{E}\|z^r-z^{r+1}\|^2+3\mathbb{E}\|u^{r+1}-z^r+\beta \nabla g(z^r)\|^2 + 3 L^2_0\beta^2\|z^{r+1}-z^r\|^2\nonumber\\
	&\stackrel{\rm(c)}\le\frac{10}{3}\mathbb{E}\|z^r-z^{r+1}\|^2+3\beta^2\mathbb{E}\left[\|\nabla g(z^r)-\frac{\lambda_{i_{r}}^r+1/N \nabla g_{i_{r}}(z^r)}{\alpha_{i_{r}}}-\sum_{i=1}^N1/N\nabla g_i(y_i^{r-1})\|^2\right]\nonumber\\
	&\stackrel{\eqref{eq:bd:var}}\leq \frac{10}{3}\mathbb{E}\|z^r-z^{r+1}\|^2+3\beta^2\sum_{i=1}^N\frac{1}{\alpha_i}\mathbb{E}\left\|\frac{1}{N}\nabla g_{i}(z^r)-\frac{1}{N}\nabla g_i(y_i^{r-1})\right\|^2\nonumber\\
	&\leq \frac{10}{3}\mathbb{E}\|z^r-z^{r+1}\|^2+3\sum_{i=1}^N\frac{\beta}{\eta_i}\mathbb{E}\left\|\frac{1}{N}\nabla g_{i}(z^r)-\frac{1}{N}\nabla g_i(y_i^{r-1})\right\|^2
	\end{align}}
where $\rm(a)$ is due to \eqref{eq:z:prox}; $\rm(b)$ is true due to the nonexpansivness of the prox operator, and the Cauchy-Swartz inequality; in $\rm(c)$ we have used the definition of $u$ in \eqref{eq:u} and the fact that $3L_0\le\sum_{i=1}^{N}\eta_i = \frac{1}{\beta} $  [cf. Assumption A-(c)]. 
In the last inequality we have applied \eqref{eq:parameter:select}, which implies that
\begin{align}\label{eq:temp}
\frac{\beta}{\alpha_i} = \frac{1}{ p_i\sum_{j=1}^{N}\eta_j} = \frac{1}{\eta_i}.
\end{align}
Note that $\eta_i$'s has to satisfy \eqref{eq:eta}. Let us follow \eqref{eq:p:eta} and choose 
$$\eta_i= \frac{9L_i}{p_i N} = \frac{9 \sum_{j=1}^{N}\eta_j}{N\eta_i} L_i.$$ 
We have
\begin{align}\label{eq:eta:equation}
\eta_i =\sqrt{9L_i/N \sum_{j=1}^{N} \eta_j}=\sqrt{9L_i/N}\sqrt{ \sum_{j=1}^{N} \eta_j}
\end{align}
Summing $i$ from $1$ to $N$ we have
\begin{align}
\sqrt{\sum_{i=1}^{N} \eta_i} = \sum_{i=1}^{N}\sqrt{9 L_i /N}
\end{align} 
Then we conclude that 
\begin{align}
\frac{1}{\beta} = {\sum_{i=1}^{N} \eta_i} = \left(\sum_{i=1}^{N}\sqrt{9 L_i /N}\right)^2.
\end{align}
So plugging the expression of $\beta$ into \eqref{eq:temp} and \eqref{eq:eta:equation}, we conclude
\begin{align}
\alpha_i = p_i = \frac{\sqrt{L_i/N}}{ \sum_{i=1}^{N}\sqrt{L_i /N}} , \quad \eta_i = \sqrt{9 L_i/N}\sum_{j=1}^{N}\sqrt{9 L_j /N}.
\end{align}

After plugging in the above inequity into \eqref{eq:gap:def}, we obtain:
{
	\begin{align}\label{eq:g:bound}
	\mathbb{E}[G^r]& \stackrel{\eqref{eq:Q:bound}}\leq \frac{10}{3\beta^2}\mathbb{E}\|z^r-z^{r+1}\|^2+\sum_{i=1}^N\frac{3}{\beta \eta_i}\mathbb{E}\left\|\frac{1}{N}\nabla g_{i}(z^r)-\frac{1}{N}\nabla g_i(y_i^{r-1})\right\|^2 \\
	&\stackrel{\eqref{eq:nesttg:descent}}\leq \frac{80}{3\beta}\mathbb{E}[Q^r-Q^{r+1} ]=\frac{80}{3} \left(\sum_{i=1}^N\sqrt{L_i/N}\right)^2\mathbb{E}[Q^r-Q^{r+1} ]\nonumber
	\end{align}}
If we sum both sides over $r=1, \cdots, R$, we obtain:
\begin{align*}
\sum_{r=1}^R\mathbb{E}[G^r]\leq \frac{80}{3} \left(\sum_{i=1}^N\sqrt{L_i/N}\right)^2\mathbb{E}[{Q^1-Q^{R+1} }].
\end{align*}
Using the definition of $z^m$, we have 
$$\mathbb{E}[G^m]=\mathbb{E}_{\cF^{r}}\left[\mathbb{E}_{m}[G^m\mid \cF^{r}] \right]=1/R \sum_{r=1}^R\mathbb{E}_{\cF^{r}}[G^r].$$
Therefore, we can finally conclude that:
\begin{align}
\mathbb{E}[G^m]\leq \frac{80}{3} \left(\sum_{i=1}^N\sqrt{L_i/N}\right)^2\frac{\mathbb{E}[{Q^1-Q^{R+1} }]}{R}\label{eq:part1}
\end{align}
which proves the first part. 

{\bf Part 2).}   In order to prove the second part let us recycle inequality in \eqref{eq:g:bound} and write
{
	\begin{align}
	&\mathbb{E}\left[G^r + \sum_{i=1}^N\frac{3}{\beta \eta_i}\left\|\frac{1}{N}\nabla g_{i}(z^{r})-\frac{1}{N}\nabla g_i(y_i^{r-1})\right\|^2\right]\nonumber\\
	&\leq \frac{10}{3\beta^2}\mathbb{E}\|z^{r+1}-z^{r}\|^2+\sum_{i=1}^N\frac{6}{\beta \eta_i}\mathbb{E}\left\|\frac{1}{N}\nabla g_{i}(z^{r})-\frac{1}{N}\nabla g_i(y_i^{r-1})\right\|^2 \nonumber\\
	&\le \frac{80}{3\beta}\mathbb{E}[Q^{r}-Q^{r+1} ]=48 \left(\sum_{i=1}^N\sqrt{L_i/N}\right)^2\mathbb{E}[Q^{r}-Q^{r+1} ]\nonumber.
	\end{align}}
Also note that 
{
	\begin{align}
	\mathbb{E}_{x^{r}}\left[\left\|x^{r+1}_i - z^{r}\right\|^2\mid \cF^r\right] = \sum_{i=1}^{N}\frac{1}{\alpha_i \eta^2_i}\left\|\frac{1}{N}\nabla g_{i}(z^{r})-\frac{1}{N}\nabla g_i(y_i^{r-1})\right\|^2
	\end{align}}
Combining the above two inequalities, we conclude
{
	\begin{align}
	&\mathbb{E}_{\cF^{r}}[G^r]+ \mathbb{E}_{\cF^{r}}\left[\sum_{i=1}^N{3 \eta^2_i }\left\|x^{r+1}_i-z^{r}\right\|^2\right]\nonumber\\
	&=\mathbb{E}_{\cF^{r}}[G^r]+ \mathbb{E}_{\cF^{r}}\left[\sum_{i=1}^N\frac{3 \eta_i \alpha_i }{\beta}\left\|x^{r+1}_i-z^{r}\right\|^2\right]\nonumber\\
	&=\mathbb{E}\left[G^r + \sum_{i=1}^N\frac{3}{\beta \eta_i}\left\|\frac{1}{N}\nabla g_{i}(z^{r})-\frac{1}{N}\nabla g_i(y_i^{r-1})\right\|^2\right]\nonumber\\
	& \le \frac{80}{3} \left(\sum_{i=1}^N\sqrt{L_i/N}\right)^2\mathbb{E}_{\cF^{r}}[Q^{r}-Q^{r+1} ]
	\end{align}}
where in the first equality we have used the relation $\frac{\alpha_i}{\beta} = \eta_i$ [cf. \eqref{eq:temp}]. 
Using a similar argument as in first part, we conclude that
\begin{align}
\mathbb{E}[G^m] + \mathbb{E}\left[\sum_{i=1}^N{3 \eta^2_i }\left\|x^{m}_i-z^{m-1}\right\|^2\right]\leq \frac{80}{3} \left(\sum_{i=1}^N\sqrt{L_i/N}\right)^2\frac{\mathbb{E}[{Q^1-Q^{R+1}} ]}{R}\label{eq:part2}.
\end{align}
This completes the proof. 
\QED

\noindent\subsection{Proof of Theorem \ref{thm:lin:conv:nesttg}}

\noindent {\bf Proof of Lemma \ref{lm:eb}}

The first statement holds true largely due to \cite[Theorem 4]{tseng09coordiate}, and the second statement holds true due to \cite[Lemma 2.1]{luo92}; see detailed discussion after \cite[Assumption 2]{tseng09coordiate}. Here the only difference with the statement  \cite[Theorem 4]{tseng09coordiate}  is that the error bound condition \eqref{eq:primaleb} holds true {\it globally}. This is by the assumption that $Z$ is a compact set. Below we provide a brief argument. 

From  \cite[Theorem 4]{tseng09coordiate}, we know that when Assumption B is satisfied, we have that for any $\xi\ge \min_z f(z)$, there exists scalars $\tau$ and $\epsilon$ such that the following error bound holds
\begin{equation}\label{eq:primaleb:2}
\dist(z,Z^*)\le \tau \|\tilde\nabla_{1/\beta} f(z)\|,\quad\mbox{whenever}~ \|\tilde\nabla_{1/\beta} f(z)\|\le \epsilon, \; f(z)\le \xi.
\end{equation}
To argue that when $Z$ is compact, the above error bound is independent of $\epsilon$,  we use the following two
steps: (1) for all $z\in Z\cap{\rm dom}(h)$ such that
$\|\tilde\nabla_{1/\beta} f(z)\|\le \delta$, it is clear that the error bound
\eqref{eq:primaleb} holds true; (2) for all $z\in
Z\cap {\rm dom}(h)$ such that $\|\tilde\nabla_{1/\beta} f(z)\|\ge \delta$,
the ratio $\frac{\dist(z,Z^*)}{ \|\tilde\nabla_{1/\beta} f(z)\|}$
is a continuous function and well defined over the compact set
$Z\cap {\rm dom}(h)\cap\left\{z\mid \|\tilde\nabla_{1/\beta} f(z)\|\ge
\delta\right\}.$ Thus, the above ratio must be bounded from  above
by a constant $\tau'$ ({independent of $b$, and no greater than ${\max_{z,z'\in Z}\|z-z'\|}/{\delta}$}). Combining (1) and (2)
yields the desired error bound over the set $Z\cap {\rm dom}(h)$. \QED

{\bf Proof of Theorem \ref{thm:lin:conv:nesttg}}

From Theorem \ref{thm:rate:nestt-g} we know that $(x^r,z^r,\lambda^r)$ converges to the set of stationary solutions of problem \eqref{eq:distributed}. Let $(x^*, z^*, \lambda^*)$ be one of such stationary solution. Then by the definition of the $Q$ function and the fact that the successive differences of the gradients goes to zero (cf. \eqref{eq:delta:go:0}), we have
\begin{align}
Q^*  = f(z^*) = \sum_{i=1}^{N}1/N g_i(z^*) + g_0(z^*)+ p(z^*).
\end{align}

Then by Lemma \ref{lm:eb} - (2) we know that $f(z^r)=\sum_{i=1}^{N}1/N g_i(z^r) + g_0(z^r)+ p(z^r)$ will finally settle at some isocost surface of $f$, i.e., there exists some {\it finite} $\bar{r}>0$ such that for all $r>\bar{r}$ and $\bar{v}\in \mathbb{R}$ such that
\begin{align}
f(\bar{z}^r) = \bar{v}, \quad \forall~r \ge \bar{r} \label{eq:vbar}
\end{align}
where $\bar{z}^r=\arg\min_{z\in Z^*}\|z^r-z\|.$ Therefore, combining the fact that $\|x^{r+1}-x^r\|\to 0$, $\|z^{r+1}-z^r\|\to 0$, $\|x^{r+1}_i - z^{r+1}\|\to 0$ and $\|\lambda^{r+1}-\lambda^r\|\to 0$ (cf. \eqref{eq:difference:z:0}, \eqref{eq:difference:lambda:0}), it is easy to see that 
\begin{align}
L(\bar{z}^r, \bar{x}^r, \bar{\lambda}^r) = f(\bar{z}^r) =\bar{v}, \quad \forall~r \ge \bar{r},
\end{align}
where $\bar{x}^r, \bar{\lambda}^r$ are defined similarly as $\bar{z}^r$. 

Now we prove that the expectation of $\Delta^{r+1} := Q^{r+1}-\bar{v}$ diminishes Q-linearly.  All the expectation below is w.r.t. the natural history of the algorithm. The proof consists of the following steps:\\
{\bf Step 1:} There exists $\sigma_1>0$ such that 
\begin{align*}
\mathbb{E}[Q^r-Q^{r+1}] \geq \sigma_1\left(\mathbb{E}\|z^{r+1}-z^r\|^2+\sum_{i=1}^N\mathbb{E}\|1/N \nabla g_i(z^r)-1/N \nabla g_i(y_i^{r-1})\|^2\right);
\end{align*}
{\bf Step 2:} There exists $\tau>0$ such that 
\begin{align*}
\mathbb{E}\|z^{r}-\bar{z}^r\|^2\leq \tau\|\mathbb{E}[\nabla_{1/\beta}\tilde{f}(z^{r})]\|^2;
\end{align*} 
{\bf Step 3:} There exists $\sigma_2>0$ such that 
\begin{align*}
\|\mathbb{E}[\nabla_{1/\beta}\tilde{f}(z^{r})]\|^2 \leq \sigma_2\left(\mathbb{E}\|z^{r+1}-z^r\|^2+\sum_{i=1}^N\mathbb{E}\|1/N \nabla g_i(z^r)-1/N \nabla g_i(y_i^{r-1})\|^2\right);
\end{align*}
{\bf Step 4:} There exists $\sigma_3>0$ such that the following relation holds true for all $r\ge \bar{r}$
\begin{align*}
\mathbb{E}[Q^{r+1}-\bar{v}] \leq \sigma_3\left(\mathbb{E}\|z^{r}-\bar{z}^r\|^2+\mathbb{E}\|z^{r+1}-z^r\|^2+\sum_{i=1}^N\mathbb{E}\|1/N \nabla g_i(z^r)-1/N \nabla g_i(y_i^{r-1})\|^2\right).
\end{align*}
These steps will be verified one by one shortly. But let us suppose that they all hold true. Below we show that linear convergence can be obtained. 

Combining step 4 and step 2 we conclude that there exists $\sigma_3>0$ such that for all $r\ge \bar{r}$
\begin{align*}
\mathbb{E}[Q^{r+1}-\bar{v}] \leq \sigma_3\left(\tau\|\mathbb{E}[\nabla_{1/\beta}\tilde{f}(z^{r-1})]\|^2+\mathbb{E}\|z^{r+1}-z^r\|^2+\sum_{i=1}^N\mathbb{E}\|1/N \nabla g_i(z^r)-1/N \nabla g_i(y_i^{r-1})\|^2\right).
\end{align*}
Then if we bound $\|\mathbb{E}(G^r)\|^2$ using step 3, we can simply make a $\sigma_4>0$ such that
\begin{align*}
\mathbb{E}[Q^{r+1}-\bar{v}] \leq \sigma_4\left(\mathbb{E}\|z^{r+1}-z^r\|^2+\sum_{i=1}^N \mathbb{E} \|1/N \nabla g_i(z^r)-1/N \nabla g_i(y_i^{r-1})\|^2\right).
\end{align*}
Finally, applying step 1 we reach the following bound for $\mathbb{E}[Q^{r+1}-\bar{v}]$:
\begin{align*}
\mathbb{E}[Q^{r+1}-\bar{v}] \leq \frac{\sigma_4}{\sigma_1}\mathbb{E}[Q^r-Q^{r+1}], \quad \forall~r\ge \bar{r},
\end{align*}
which further implies that for $\sigma_5 = \frac{\sigma_4}{\sigma_1}>0$, we have
\begin{align*}
\mathbb{E}[\Delta^{r+1}]\leq \frac{\sigma_5}{1+\sigma_5}\mathbb{E}[\Delta^r], \quad \forall~r\ge \bar{r}. 
\end{align*} 

Now let us verify the correctness of each step. Step 1 can be directly obtained from equation \eqref{eq:nesttg:descent}. Step 2 is exactly Lemma \eqref{lm:eb}. Step 3 can be verified using a similar derivation as in \eqref{eq:Q:bound}\footnote{We simply need to  replace $-z^{r-1}+\prox_h^{1/\beta}[u^{r-1} -\beta \nabla g_0(z^{r-1})]$ in step (a) of \eqref{eq:Q:bound} by $-z^{r}+\prox_h^{1/\beta}[u^{r} -\beta \nabla g_0(z^{r})]$ and using the same derivation.}. 

Below let us prove the step 4, which is a bit involved. From \eqref{eq:z:nestt} we know that 
\begin{align*}
z^{r+1}= \arg\min_z h(z)+g_0(z)+\sum_{i=1}^N \langle \lambda_i^{r},x^{r+1}_i-z\rangle +\frac{\eta_i}{2}\|x_i^{r+1}-z\|^2.
\end{align*}
This implies that
\begin{align}
&h(z^{r+1})+g_0(z^{r+1})+\sum_{i=1}^N \langle \lambda_i^{r},x^{r+1}_i-z^{r+1}\rangle +\frac{\eta_i}{2}\|x_i^{r+1}-z^{r+1}\|^2\nonumber\\
&\le h(\bar{z}^r)+g_0(\bar{z}^r)+\sum_{i=1}^N \langle \lambda_i^{r},x^{r+1}_i-\bar{z}^r\rangle +\frac{\eta_i}{2}\|x_i^{r+1}-\bar{z}^r\|^2.
\end{align}
Rearranging the terms, we obtain
\begin{align*}
h(z^{r+1})+g_0(z^{r+1})-h(\bar{z}^r)-g_0(\bar{z}^r)\leq \sum_{i=1}^N \langle \lambda_i^{r},z^{r+1}-\bar{z}^r\rangle +\frac{\eta_i}{2}\|x_i^{r+1}-\bar{z}^r\|^2.
\end{align*}
Using this inequality we have:
\begin{align}\label{eq:QV}
Q^{r+1}-\bar{v}&\leq \sum_{i=1}^N 1/N\left(g_i(z^{r+1})-g_i(\bar{z}^r)\right) + \langle \lambda_i^{r},z^{r+1}-\bar{z}^r\rangle\nonumber \\
&+\sum_{i=1}^N\frac{\eta_i}{2}\|x_i^{r+1}-\bar{z}^r\|^2 + \|1/N(\nabla g_i(z^{r})-\nabla g_i(y_i^{r-1})\|^2.
\end{align}
The first term in RHS can be bounded as follows:
\begin{align*}
&\sum_{i=1}^N 1/N\left(g_i(z^{r+1})-g_i(\bar{z}^r)\right)\nonumber\\
&\stackrel{(a)}\leq \sum_{i=1}^N 1/N \langle \nabla g_i(\bar{z}^r), z^{r+1}-\bar{z}^r\rangle +L_i/2N \|z^{r+1}-\bar{z}^r\|^2\\
&\leq \sum_{i=1}^N 1/N \langle \nabla g_i(\bar{z}^r)+\nabla g_i(z^{r+1})-\nabla g_i(z^{r+1}), z^{r+1}-\bar{z}^r\rangle +L_i/2N \|z^{r+1}-\bar{z}^r\|^2\\
&\stackrel{(b)}\leq \sum_{i=1}^N 1/N \langle \nabla g_i(z^{r+1}), z^{r+1}-\bar{z}^r\rangle + 3L_i/2N \|z^{r+1}-\bar{z}^r\|^2,
\end{align*}
where $(a)$ is true due to the descent lemma; and $(b)$ comes from the Lipschitz continuity of the $\nabla g_i$. 

Plugging the above bound into \eqref{eq:QV}, we further have:
\begin{align}
Q^{r+1}-\bar{v}&\leq \sum_{i=1}^N 1/N\langle \nabla g_i(z^{r+1})-\nabla g_i(y_i^{r-1}), z^{r+1}-\bar{z}^r\rangle + 3L_i/2N \|z^{r+1}-\bar{z}^r\|^2\nonumber\\
&+\frac{\eta_i}{2}\|x_i^{r+1}-\bar{z}^r\|^2 + \|1/N(\nabla g_i(z^{r})-\nabla g_i(y_i^{r-1})\|^2\nonumber\\
&=\sum_{i=1}^N 1/N\langle \nabla g_i(z^{r+1})+\nabla g_i(z^{r})-\nabla g_i(z^{r})-\nabla g_i(y_i^{r-1}), z^{r+1}-\bar{z}^r\rangle\nonumber\\
&+\frac{\eta_i}{2}\|x_i^{r+1}-\bar{z}^r\|^2 + \|1/N(\nabla g_i(z^{r})-\nabla g_i(y_i^{r-1})\|^2+  3L_i/2N \|z^{r+1}-\bar{z}^r\|^2,\nonumber
\end{align}
where in the first inequality we have used the fact that $\lambda^r_i = -\frac{1}{N}\nabla g_i(y^{r-1}_i)$; cf . \eqref{eq:lambda:compact}.
Applying the Cauchy-Schwartz inequality we further have:
\begin{align}
Q^{r+1}-\bar{v}&\leq\sum_{i=1}^N 1/2\|1/N\left(\nabla g_i(z^{r+1})+\nabla g_i(z^{r})\right)\|^2+1/2\|z^{r+1}-\bar{z}^r\|^2\nonumber\\
&+ \sum_{i=1}^N 1/2\|1/N\left(\nabla g_i(z^{r})-\nabla g_i(y_i^{r-1})\right)\|^2+1/2\|z^{r+1}-\bar{z}^r\|^2\nonumber\\
&+\frac{\eta_i}{2}\|x_i^{r+1}-\bar{z}^r\|^2 + \|1/N(\nabla g_i(z^{r})-\nabla g_i(y_i^{r-1})\|^2 + 3L_i/2N \|z^{r+1}-\bar{z}^r\|^2\nonumber\\
&\leq \sum_{i=1}^N\left[ \frac{L_i^2}{2N^2}\|z^{r+1}-z^{r}\|^2+\frac{3}{2N^2}\|g_i(z^{r})-\nabla g_i(y_i^{r-1})\|^2+\frac{\eta_i}{2}\|x_i^{r+1}-\bar{z}^r\|^2\right]\nonumber\\
&+\left(1+ 3L_i/2N \right)\|z^{r+1}-\bar{z}^r\|^2  \label{eq:dif:Q}.
\end{align}
Now let us bound $\sum_{i=1}^N\frac{\eta_i}{2}\|x_i^{r+1}-\bar{z}^r\|^2$ in the above inequality:
\begin{align}
\sum_{i=1}^N\frac{\eta_i}{2}\|x_i^{r+1}-\bar{z}^r\|^2&=\sum_{i=1}^N\frac{\eta_i}{2}\|x_i^{r+1}-z^{r+1}+z^{r+1}-\bar{z}^r\|^2\nonumber\\
&\leq \sum_{i=1}^N\eta_i\|x_i^{r+1}-z^{r+1}\|^2+\eta_i\|z^{r+1}-\bar{z}^r\|^2\nonumber\\
&=\sum_{i=1}^N\eta_i\|x_i^{r+1}-z^{r}+z^{r}-z^{r+1}\|^2+\eta_i\|z^{r+1}-\bar{z}^r\|^2\nonumber\\
&\leq \sum_{i=1}^N2\eta_i\|x_i^{r+1}-z^{r}\|^2+2\eta_i\|z^{r}-z^{r+1}\|^2+\eta_i\|z^{r+1}-\bar{z}^r\|^2\nonumber.
\end{align}
Using the fact that $x_i^{r+1}=z^r$ when $i\neq i_{r}$ we further have:
\begin{align}
\sum_{i=1}^N\frac{\eta_i}{2}\|x_i^{r+1}-\bar{z}^r\|^2&\leq 2\eta_{i_{r}}\|x_{i_{r}}^{r+1}-z^{r}\|^2+\sum_{i=1}^N2\eta_i\|z^{r}-z^{r+1}\|^2+\eta_i\|z^{r+1}-\bar{z}^r\|^2\nonumber\\
&= \frac{2}{\alpha_{i_{r}}^2\eta_{i_{r}}}\|\lambda_{i_{r}}+1/N\nabla g_{i_{r}}(z^{r})\|^2+\sum_{i=1}^N2\eta_i\|z^{r}-z^{r+1}\|^2+\eta_i\|z^{r+1}-\bar{z}^r\|^2\nonumber\\
&= \frac{2}{\alpha_{i_{r}}^2\eta_{i_{r}}N^2}\|\nabla g_{i_{r}}(z^{r})-\nabla g_{i_{r}}(y_{i_{r}}^{r-1})\|^2\nonumber\\
&+\sum_{i=1}^N2\eta_i\|z^{r}-z^{r+1}\|^2+\eta_i\|z^{r+1}-z^r+z^r-\bar{z}^r\|^2\nonumber\\
&\leq \frac{2}{\alpha_{i_{r}}^2\eta_{i_{r}}N^2}\|\nabla g_{i_{r}}(z^{r})-\nabla g_{i_{r}}(y_{i_{r}}^{r-1})\|^2\nonumber\\
&+\sum_{i=1}^N4\eta_i\|z^{r}-z^{r+1}\|^2+2\eta_i\|z^r-\bar{z}^r\|^2\label{eq:term1}.
\end{align}
Take expectation on both sides of the above equation and set $p_i=\alpha_i$, we obtain:
\begin{align*}
\sum_{i=1}^N\frac{\eta_i}{2}\mathbb{E}\|x_i^{r+1}-\bar{z}^r\|^2
&\leq \sum_{i=1}^N\frac{2}{\alpha_{i}\eta_{i}}\mathbb{E}\|\nabla g_{i}(z^{r})-\nabla g_{i}(y_{i}^{r-1})\|^2\\
&+\sum_{i=1}^N4\eta_i\mathbb{E}\|z^{r}-z^{r+1}\|^2+2\eta_i\mathbb{E}\|z^r-\bar{z}^r\|^2.
\end{align*}
Combining equations \eqref{eq:dif:Q} and \eqref{eq:term1}, eventually one can find $\sigma_3>0$ such that 
\begin{align*}
\mathbb{E}[Q^{r+1}-\bar{v}] \leq \sigma_3\left(\mathbb{E}\|z^{r}-\bar{z}\|^2+\mathbb{E}\|z^{r+1}-z^r\|^2+\sum_{i=1}^N\mathbb{E}\|1/N \nabla g_i(z^r)-1/N \nabla g_i(y_i^{r-1})\|^2\right),
\end{align*}
which completes the proof of Step 4.

In summary, we have shown that Step 1 - 4 all hold true. Therefore we have shown that the NESTT-G converges Q-linearly.   
\QED

\noindent\subsection{Some Key Properties of NESTT-E}
To facilitate the following derivation, in this section we collect some key properties of NESTT-E. 

First, for $i=i_{r}$, using the optimality condition for $x_i$ update step \eqref{eq:xi_update:exact} we have the following identity:
\begin{align}\label{eq:nestte:x:opt}
\frac{1}{N}\nabla g_{i_{r}}(x_{i_{r}}^{r+1})+\lambda^r_{i_{r}}+\alpha_{i_{r}}\eta_{i_{r}} (x^{r+1}_{i_{r}}-z^{r+1})=0.
\end{align}
Combined with the dual variable update step \eqref{eq:y_update:exact} we obtain
\begin{align}
\frac{1}{N}\nabla g_{i_{r}}(x_{i_{r}}^{r+1})=-\lambda^{r+1}_{i_{r}}\label{eq:y_exp2}.
\end{align}

Second, the optimality condition for the $z$-update is given by: 
\begin{align}
z^{r+1}&=\prox_h\left[ z^{r+1}-\nabla_z (L(x^r,z,\lambda^r)-h(z))\right]\\
&=\prox_h\left[ z^{r+1}-\sum_{i=1}^{N}\eta_i\left(z^{r+1}-x^r_i-\frac{\lambda^r_i}{\eta_i}\right)  - \nabla g_0(z^{r+1}) \right] \label{eq:optcond:z}.
\end{align}

\noindent\subsection{Proof of Theorem \ref{thm:sublin}}
To prove this result, we need a few lemmas. 

For notational simplicity, define new variables $\{\hx^{r+1}_i\}$,  $\{\hat{\lambda}_i^{r+1}\}$ by
\begin{align}
\hx_i^{r+1}&:=\arg\min_{x_i} \; U_i(x_i,z^{r+1},\lambda^r_i), \quad \hat{\lambda}^{r+1}_i:=\lambda_{i}^{r}+\alpha_i\eta_i\left(\hx^{r+1}_i-z^{r+1}\right), \quad \forall i.\label{eq:hy}
\end{align}
These variables are the {\it virtual variables} generated by updating all variables at iteration $r+1$. Also define:
\begin{align*}
&L^r:= L(x^r, z^{r}; \lambda^{r}), \quad w := (x, z, \lambda), \quad \beta:=\frac{1}{\sum_{i=1}^N\eta_i}, \quad c_i:= \frac{L^2_i}{\alpha_i\eta_i N^2}-\frac{\gamma_i}{2}+\frac{1-\alpha_i}{\alpha_i}\frac{L_i}{N}
\end{align*}

First, we need the following lemma to show that the size  of the successive difference of the dual variables can be upper bounded by that of the primal variables. This is a simple consequence of \eqref{eq:y_exp2}; also see [R2, Lemma 2.1]. We include the proof for completeness. 
\begin{lemma}\label{lemma:y1}
	Suppose assumption A holds. Then for NESTT-E algorithm, the following are true:
	\begin{subequations}
		\begin{align}
		&\|\lambda^{r+1}_i-\lambda_{i}^{r}\|^2 \leq \frac{L^2_i}{N^2}\|x_i^{r+1}-x_{i}^{r}\|^2, \quad \|\hat{\lambda}^{r+1}_i-\lambda_{i}^{r}\|^2 \leq \frac{L^2_i}{N^2}\|\hx_i^{r+1}-x_{i}^{r}\|^2,\; \forall~i.\label{eq:y_difference_rand1}
		\end{align}
	\end{subequations}
\end{lemma}

{\bf Proof.}  {We only show the first inequality. The second one follows an analogous argument.
	
	To prove \eqref{eq:y_difference_rand1}, first note that the case for $i\neq i_{r}$ is trivial, as both sides of \eqref{eq:y_difference_rand1} are zero. For the index $i_{r}$, we have a closed-form expression for $\lambda^r_{i_{r}}$ following \eqref{eq:y_exp2}. 
	Notice  that for any given $i$, the primal-dual pair $(x_{i}, \lambda_i)$ is always updated at the same iteration. Therefore,  if for each $i$ we choose the initial solutions in a way such that $\lambda_i^0=-\nabla g_i(x_i^0)$, then we have  
	\begin{align}\label{eq:y_expression1}
	\frac{1}{N}\nabla g_i(x_i^{r+1})=-\lambda^{r+1}_i \quad \forall~ i=1,2,\cdots N.
	\end{align}
	
	Combining \eqref{eq:y_expression1} with Assumption A-(a) yields the following:
	\begin{align*}
	\|\lambda^{r+1}_i-\lambda^r_i\|= \frac{1}{N}\|\nabla g_i(x_i^{r+1})-\nabla g_i(x_i^{r})\|\leq \frac{L_i}{N}\|x_i^{r+1}-x^{r}_i\|.
	\end{align*}
	The proof is complete.}\QED

Second, we bound the successive difference of the potential function. 
\begin{lemma}\label{lemma:L_diff:exact}
	Suppose Assumption A holds true. Then the following holds for NESTT-E
	\begin{align}\label{eq:L_descent1:exact}
	&\mathbb{E}[L^{r+1}-L^r|x^r, z^r]\leq -\frac{\gamma_z}{2}\|z^{r+1}-z^r\|^2 +\sum_{i=1}^Np_ic_{i}\|x_{i}^r-\hat{x}_{i}^{r+1}\|^2.
	\end{align}
\end{lemma}

{\bf Proof.}   First let us split $L^{r+1}-L^r$ in the following way:
\begin{align}\label{eq:successive_L1:exact}
&L^{r+1}-L^r=L^{r+1}-L(x^{r+1}, z^{r+1}; \lambda^{r})+L(x^{r+1}, z^{r+1}; \lambda^{r})-L^r.
\end{align}
The first two terms in \eqref{eq:successive_L1:exact} can be bounded by
\begin{align}
&L^{r+1}-L(x^{r+1}, z^{r+1}; \lambda^{r})=\sum_{i=1}^{N}\langle \lambda_i^{r+1}-\lambda_{i}^{r}, x^{r+1}_i-z^{r+1}\rangle\nonumber\\
&\stackrel{\rm (a)}=\frac{1}{\alpha_{i_{r}}\eta_{i_{r}}}\|\lambda_{i_{r}}^{r+1}-\lambda_{i_{r}}^{r}\|^2\stackrel{\rm (b)}\leq \frac{L^2_{i_{r}}}{N^2\alpha_{i_{r}}\eta_{i_{r}}}\|x^{r+1}_{i_{r}}-x_{{i_{r}}}^{r}\|^2\label{eq:extra:exact}
\end{align}
where in $\rm (a)$ we have used \eqref{eq:y_update:exact}, and the fact that $\lambda_i^{r+1}-\lambda_{i}^{r}=0$ for all variable blocks except $i_{r}$th block; $\rm (b)$ is true because of Lemma \ref{lemma:y1}. 

The last two terms in \eqref{eq:successive_L1:exact} can be written in the following way:
\begin{align}
&L(\{x^{r+1}_i\}, z^{r+1}; \lambda^{r})-L^r\nonumber\\
&=L(x^{r+1}, z^{r+1}; \lambda^{r})-L(x^{r}, z^{r+1}; \lambda^{r})+L(x^{r}, z^{r+1}; \lambda^{r})-L^r\label{eq:L:z:diff}.
\end{align}

The first two terms in \eqref{eq:L:z:diff} characterizes the change of the Augmented Lagrangian before and after the update of $x$. Note that $x$ updates do not directly optimize the augmented Lagrangian. Therefore the characterization of this step is a bit involved. We have the following:
\begin{align}
&L(x^{r+1}, z^{r+1}; \lambda^{r})-L(x^{r}, z^{r+1}; \lambda^{r})\nonumber\\
&\stackrel{\rm (a)}\le \sum_{i=1}^{N}\left(\left\langle\nabla_{i}  L(x^{r+1}, z^{r+1}; \lambda^{r}), x^{r+1}_i-x^{r}_i\right\rangle-\frac{{\gamma}_i}{2}\|x^{r+1}_i-x_{i}^{r}\|^2\right)\nonumber\\
&\stackrel{\rm (b)}= \left\langle\nabla_{i_{r}}  L(x^{r+1}, z^{r+1}; \lambda^{r}), x^{r+1}_{i_{r}}-x^{r}_{i_{r}}\right\rangle-\frac{{\gamma}_{i_{r}}}{2}\|x^{r+1}_{i_{r}}-x_{{i_{r}}}^{r}\|^2\nonumber\\
&\stackrel{\rm (c)}= \left\langle\eta_{i_{r}}(1-\alpha_{i_{r}})(x^{r+1}_{{i_{r}}}-z^{r+1}),x^{r+1}_{{i_{r}}}-x^{r}_{{i_{r}}}\right\rangle -\frac{{\gamma}_{i_{r}}}{2}\|x^{r+1}_{i_{r}}-x_{{i_{r}}}^{r}\|^2\nonumber\\
&\stackrel{\rm (d)}= \left\langle\frac{1-\alpha_{i_{r}}}{\alpha_{i_{r}}}(\lambda^{r+1}_{{i_{r}}}-\lambda^{r}_{i_{r}}),x^{r+1}_{{i_{r}}}-x^{r}_{{i_{r}}}\right\rangle -\frac{{\gamma}_{i_{r}}}{2}\|x^{r+1}_{i_{r}}-x_{{i_{r}}}^{r}\|^2\nonumber\\
&\le \frac{1-\alpha_{i_{r}}}{\alpha_{i_{r}}} \left(\frac{1}{2 L_{i_{r}}/N} \|\lambda^{r+1}_{{i_{r}}}-\lambda^{r}_{i_{r}}\|^2+ \frac{L_{i_{r}}}{2 N}\|x^{r+1}_{{i_{r}}}-x^{r}_{{i_{r}}}\|^2\right)-\frac{{\gamma}_{i_{r}}}{2}\|x^{r+1}_{i_{r}}-x_{{i_{r}}}^{r}\|^2\nonumber\\
&\stackrel{\rm (e)}\le \frac{1-\alpha_{i_{r}}}{\alpha_{i_{r}}}\frac{L_{i_{r}}}{N}\|x^{r+1}_{{i_{r}}}-x^{r}_{{i_{r}}}\|^2 -\frac{{\gamma}_{i_{r}}}{2}\|x^{r+1}_{i_{r}}-x_{{i_{r}}}^{r}\|^2\label{eq:L:x_i:diff2}
\end{align}
where
\begin{itemize}
	\item $\rm(a)$ is true because $L(x,z,\lambda)$ is strongly convex with respect to $x_i$.
	\item $\rm(b)$ is true because when $i\neq i_{r}$, we have $x^{r+1}_i=x_{i}^{r}$.
	\item $\rm(c)$ is true because $x^{r+1}_{i_{r}}$ is optimal solution for the problem $\min U_{i_{r}}(x_{i_{r}},z^{r+1},\lambda_{i_{r}}^{r})$ (satisfying \eqref{eq:nestte:x:opt}), and we have used the optimality of such $x^{r+1}_{i_{r}}$.
	\item $\rm(d)$ and $\rm(e)$ are due to Lemma \ref{lemma:y1}.
\end{itemize}

{Similarly, the last two terms in \eqref{eq:L:z:diff} can be  bounded using equation \eqref{eq:nestte:x:opt} and the strong convexity of function $L$ with respect to the variable $z$. Therefor We have:}
\begin{align}
L(x^r,z^{r+1},\lambda^r)-L^r\leq -\frac{\gamma_z}{2}\|z^{r+1}-z^r\|^2 \label{eq:L:z:diff2}.
\end{align}

Combining equations \eqref{eq:extra:exact}, \eqref{eq:L:x_i:diff2} and \eqref{eq:L:z:diff2}, eventually we have:
\begin{align}
&L^{r+1}-L(x^r,z^{r+1},\lambda^r)\le c_{i_{r}}\|x_{i_{r}}^r-x_{i_{r}}^{r+1}\|^2 \label{eq:key:2}\\
&L^{r+1}-L^r\le -\frac{\gamma_z}{2}\|z^{r+1}-z^r\|^2 +c_{i_{r}}\|x_{i_{r}}^r-x_{i_{r}}^{r+1}\|^2
\end{align}
Taking expectation on both side of this inequality with respect to $i_{r}$, we can conclude that:
\begin{align}
&\mathbb{E}[L^{r+1}-L^r\mid z^r, x^r]\leq -\frac{\gamma_z}{2}\|z^{r+1}-z^r\|^2 +\sum_{i=1}^Np_ic_{i}\|x_{i}^r-\hat{x}_{i}^{r+1}\|^2
\end{align}
where $p_i$ is the probability of picking $i$th block. The lemma is proved. \QED

\begin{lemma}\label{lem:pot:bd:below}
	Suppose that Assumption A is satisfied, then $L^r\geq \underline{f}$.  
\end{lemma}
{\bf Proof.} Using the definition of the augmented Lagrangian function we have:
\begin{align}
&L^{r+1}=\sum_{i=1}^{N}\left( \frac{1}{N}g_i(x^{r+1}_i)+\langle \lambda^{r+1}_i, x^{r+1}_i-z^{r+1}\rangle+\frac{\eta_i}{2}\|x^{r+1}_i-z^{r+1}\|^2\right)+g_0(z^{r+1})+p(z^{r+1})\nonumber\\
&\stackrel{\rm (a)}=\sum_{i=1}^{N}\left( \frac{1}{N} g_i(x^{r+1}_i)+\frac{1}{N} \langle \nabla g_i(x^{r+1}_i), z^{r+1}-x_i^{r+1}\rangle+\frac{\eta_i}{2}\|x^{r+1}_i-z^{r+1}\|^2\right)+g_0(z^{r+1})+p(z^{r+1})\nonumber\\
&\stackrel{\rm (b)}\ge \sum_{i=1}^{N} \frac{1}{N}g_i(z^{r+1}) +\left(\frac{\eta_i}{2}-\frac{L_i}{2 N}\right)\|z^{r+1}-x_i^{r+1}\|^2+ g_0(z^{r+1})+p(z^{r+1})\nonumber\\
&\stackrel{\rm (c)} \ge \sum_{i=1}^{N} \frac{1}{N} g_i(z^{r+1})+g_0(z^{r+1})+p(z^{r+1})\geq \underline{f} \label{eq:lb:aug}
\end{align}
where $\rm(a)$ is true because of equation \eqref{eq:y_exp2}; $\rm (b)$ follows Assumption A-(b); $\rm(c)$ follows Assumption A-(d).
The desired result is proven.\QED

{\bf Proof of Theorem \ref{thm:sublin}.} We first show that the algorithm converges to the set of stationary solutions, and then establish the convergence rate. 

{\bf Step 1. Convergence to Stationary Solutions}.  Combining the descent estimate in Lemma \ref{lemma:L_diff:exact} as well as the lower bounded condition in Lemma \ref{lem:pot:bd:below}, we can again apply the Supermartigale Convergence Theorem \cite[Proposition 4.2]{bertsekas97} and conclude that
\begin{align}\label{eq:difference:z:0}
\|x_i^{r+1}-x_i^r\| \to 0, \quad \|z^{r+1}-z^r\|\to 0, \mbox{with probability 1}.
\end{align}
From Lemma \ref{lemma:y1} we have that the constraint violation is satisfied
\begin{align}\label{eq:difference:lambda:0}
\|\lambda^{r+1}-\lambda^r\|\to 0, \quad \|x^{r+1}_i-z^r\|\to 0. 
\end{align}
The rest of the proof follows similar lines as in \cite[Theorem 2.4]{hong14nonconvex_admm}. Due to space limitations we omit the proof. 

{\bf Step 2. Convergence Rate.} We first show that there exists a $\sigma_1(\alpha)>0$ such that
\begin{align}
\|\tilde{\nabla}L(w^r)\|^2+\sum_{i=1}^{N}\frac{L_i^2}{N^2}\|x^r_i-z^r\|^2 \leq \sigma_1(\alpha)\left(\|z^r-z^{r+1}\|^2+\sum_{i=1}^{N}\|x^r_i-\hat{x}^{r+1}_i\|^2\right). \label{eq:total:square}
\end{align}
Using the definition of $\|\tilde{\nabla}L^r(w^r)\|$ we have:
\begin{align}
\|\tilde{\nabla}L^r(w^r)\|^2 &= \|z^r-\prox_{h}\left[z^r-\nabla_{z}(L^r-h(z^r))\right]\|^2\nonumber \\
&\quad+ \sum_{i=1}^N\left\|\frac{1}{N}\nabla g_i(x_i^r)+\lambda_i^r+\eta_i(x_i^r-z^r)\right\|^2 \label{eq:gap:z}.
\end{align}

{From the optimality condition of the $z$ update \eqref{eq:optcond:z} we have:}
$$z^{r+1}=\prox_h\left[ z^{r+1}-\sum_{i=1}^{N}\eta_i\left(z^{r+1}-x^r_i-\frac{\lambda^r_i}{\eta_i}\right)  - \nabla g_0(z^{r+1}) \right].$$
Using this, the first term in equation \eqref{eq:gap:z} can be bounded as: 
\begin{align}
&\left\|z^r-\prox_{h}\left[z^r-\nabla_{z}(L^r-h(z^r))\right]\right\|\nonumber\\
&=\left\|z^{r}-z^{r+1}+z^{r+1}-\prox_h\left[z^r-\sum_{i=1}^{N}\eta_i(z^r-x^{r}_i-\frac{\lambda^r_i}{\eta_i})- \nabla g_0(z^{r})\right]\right\|\nonumber\\
&\le \|z^{r}-z^{r+1}\|+\Bigg\|\prox_h\left[ z^{r+1}-\sum_{i=1}^{N}\eta_i\left(z^{r+1}-x^r_i-\frac{\lambda^r_i}{\eta_i}\right) - \nabla g_0(z^{r+1}) \right] \nonumber\\
&\quad \quad -\prox_h\left[z^r-\sum_{i=1}^{N}\eta_i(z^r-x^{r}_i-\frac{\lambda^r_i}{\eta_i})- \nabla g_0(z^{r})\right]\Bigg\|\nonumber\\
&\le 2\|z^{r+1}-z^r\|+\left(\sum_{i=1}^{N}\eta_i+L_0\right)\|z^r-z^{r+1}\|, \label{eq:prox1}
\end{align}
where in the last inequality we have used the nonexpansiveness of the proximity operator.

Similarly, the optimality condition of the $x_i$ subproblem is given by
\begin{align}
\frac{1}{N}\nabla g_i(\hat{x}_i^{r+1})+\lambda_i^r+\alpha_i\eta_i(\hat{x}_i^{r+1}-z^{r+1})=0 \label{eq:opt:x}.
\end{align} 
Applying this identity, the second term in equation \eqref{eq:gap:z} can be written as follows:
\begin{align}
&\sum_{i=1}^N\left\|\frac{1}{N}\nabla g_i(x_i^r)+\lambda_i^r+\eta_i(x_i^r-z^r)\right\|^2\nonumber\\
&\stackrel{\rm(a)}=\sum_{i=1}^N\left\|\frac{1}{N}\nabla g_i(x_i^r)-\frac{1}{N}\nabla g_i(\hat{x}_i^{r+1})+\eta_i(x_i^r-z^r)-\alpha_i\eta_i(\hat{x}_i^{r+1}-z^{r+1})\right\|^2\nonumber\\
&=\sum_{i=1}^N\left\|\frac{1}{N}\nabla g_i(x_i^r)-\frac{1}{N}\nabla g_i(\hat{x}_i^{r+1})+\eta_i(x_i^r-\hat{x}_i^{r+1}+\hat{x}_i^{r+1}-z^{r+1}+z^{r+1}-z^r)-\alpha_i\eta_i(\hat{x}_i^{r+1}-z^{r+1})\right\|^2\nonumber\\
&\stackrel{\rm(b)}\leq 4\sum_{i=1}^N\left[\left(\frac{L_i^2}{N^2}+\eta_i^2+\frac{(1-\alpha_i)^2L_i^2}{N^2\alpha_i^2}\right)\|\hat{x}_i^{r+1}-x_i^r\|^2+\eta_i^2\|z^{r+1}-z^r\|^2\right]\label{eq:prox2},
\end{align}
where $\rm(a)$ holds because of equation \eqref{eq:opt:x}; $\rm(b)$ holds because of Lemma \ref{lemma:y1}.  

Finally, combining \eqref{eq:prox1} and \eqref{eq:prox2} leads to the following bound for proximal gradient
\begin{align}
\|\tilde{\nabla}L^r\|^2&\le \left(4\sum_{i=1}^N\eta_i^2+\left({2}+L_0+\sum_{i=1}^{N}\eta_i\right)^2\right)\|z^r-z^{r+1}\|^2\nonumber\\
&+\sum_{i=1}^{N}4\left(\frac{L_i^2}{N^2}+\eta_i^2+\frac{(1-\alpha_i)^2L_i}{N^2\alpha_i^2}\right)\|x^r_i-\hat{x}^{r+1}_i\|^2.\label{eq:estimate:nL}
\end{align}

Also note that:
\begin{align}
\sum_{i=1}^N \frac{L_i^2}{N^2}\|x_i^r-z^r\|^2&\leq \sum_{i=1}^N 3\frac{L_i^2}{N^2}\left[\|x^r_i-\hat{x}^{r+1}_i\|^2+\|\hat{x}_i^{r+1}-z^{r+1}\|^2+\|z^{r+1}-z^r\|^2\right]\nonumber\\
&=\sum_{i=1}^N 3\frac{L_i^2}{N^2}\left[\|x^r_i-\hat{x}^{r+1}_i\|^2+\frac{1}{\alpha_i^2\eta_i^2}\|\hat{\lambda}_i^{r+1}-\lambda_i^{r}\|^2+\|z^{r+1}-z^r\|^2\right]\nonumber\\
&\leq\sum_{i=1}^N 3\frac{L_i^2}{N^2}\left[\|x^r_i-\hat{x}^{r+1}_i\|^2+\frac{L_i^2}{\alpha_i^2\eta_i^2 N^2}\|\hat{x}_i^{r+1}-x_i^{r}\|^2+\|z^{r+1}-z^r\|^2\right].
\label{eq:estimate:x}
\end{align}
The two inequalities \eqref{eq:estimate:nL} -- \eqref{eq:estimate:x} imply that:
\begin{align}
& \|\tilde{\nabla}L^r\|^2+\sum_{i=1}^{N}\frac{L_i^2}{N^2}\|x^r_i-z^r\|^2 \nonumber\\
&\leq \left(\sum_{i=1}^N4\eta_i^2+ ({2}+\sum_{i=1}^N\eta_i+ L_0)^2+3\sum_{i=1}^N\frac{L_i^2}{N^2}\right)\|z^r-z^{r+1}\|^2\nonumber\\
&+\sum_{i=1}^{N}\left(4\left(\frac{L_i^2}{N^2}+\eta_i^2+(\frac{1}{\alpha_i}-1)^2\frac{L^2_i}{N^2}\right)+3\left(\frac{L_i^4}{\alpha_i N^4 \eta_i^2}+\frac{L_i^2}{N^2}\right)\right)\|x^r_i-\hat{x}^{r+1}_i\|^2.
\end{align}
Define the following quantities:
\begin{align*}
\hat{\sigma}_1(\alpha) &= \max_i\left\{4\left(\frac{L_i^2}{N^2}+\eta_i^2+\left(\frac{1}{\alpha_i}-1\right)^2\frac{L^2_i}{N^2}\right)+3\left(\frac{L_i^4}{\alpha_i \eta_i^2 N^4}+\frac{L_i^2}{N^2}\right)\right\}\\
\tilde{\sigma}_1&=\sum_{i=1}^N4\eta_i^2+({2}+\sum_{i=1}^N\eta_i+ L_0 )^2+3\sum_{i=1}^N\frac{L_i^2}{N^2}.
\end{align*}
Setting  $\sigma_1(\alpha)=\max(\hat{\sigma}_1(\alpha), \tilde{\sigma}_1)>0$, we have
\begin{align}
\|\tilde{\nabla}L^r\|^2+\sum_{i=1}^{N}\frac{L_i^2}{N^2}\|x^r_i-z^r\|^2 \leq \sigma_1(\alpha)\left(\|z^r-z^{r+1}\|^2+\sum_{i=1}^{N}\|x^r_i-\hat{x}^{r+1}_i\|^2\right).
\end{align}
From Lemma \ref{lemma:L_diff:exact} we know that
\begin{align}
&\mathbb{E}[L^{r+1}-L^r|z^r, x^r]\leq -\frac{\gamma_z}{2}\|z^{r+1}-z^r\|^2 +\sum_{i=1}^Np_ic_{i}\|x_{i}^r-\hat{x}_{i}^{r+1}\|^2
\end{align}
Note that $\gamma_z=\sum_{i=1}^N \eta_i - L_0$, then define  $\hat{\sigma}_2$ and $\tilde{\sigma}_2$ as 
\begin{align*}
\hat{\sigma}_2(\alpha) &= \max_i \left\{p_i\left(\frac{\gamma_i}{2}-\frac{L^2_i}{\alpha_i\eta_i N^2}-\frac{1-\alpha_i}{\alpha_i}\frac{L_i}{N}\right)\right\}\\
\tilde{\sigma}_2 &= \frac{\sum_{i=1}^N \eta_i - L_0}{2}.
\end{align*}
We can set $\sigma_2(\alpha)=\max(\hat{\sigma}_2(\alpha), \tilde{\sigma}_2)$ to obtain
\begin{align}
&E[L^r-L^{r+1}|x^r, z^r]\ge\sigma_2(\alpha)\left(\sum_{i=1}^{N}\|\hat{x}^{r+1}_i-x_{i}^{r}\|^2+\|z^{r+1}-z^r\|^2\right).\label{eq:estimate:L}
\end{align}
{ Combining \eqref{eq:total:square} and \eqref{eq:estimate:L} we have
	\begin{align}
	H(w^r)=\|\tilde{\nabla}L^r\|^2+\sum_{i=1}^{N}L_i^2/N\|x^r_i-z^r\|^2 \leq\frac{\sigma_1(\alpha)}{\sigma_2(\alpha)}E[L^r-L^{r+1}|F^r]\nonumber.
	\end{align}}
Let us set $C(\alpha) = \frac{\sigma_1(\alpha)}{\sigma_2(\alpha)}$ and take expectation on both side of the above equation  to obtain:
\begin{align}
\mathbb{E}[H(w^r)] \leq C(\alpha) E[L^r-L^{r+1}].
\end{align}

Summing both sides of the above inequality over $r=1,\cdots, R$, we obtain:
\begin{align}
&\sum_{r=1}^{R}\mathbb{E}[H(w^r)]\leq C(\alpha)E[L^1-L^{R+1}].
\end{align}
Using the definition of $w^m=(x^m,z^m,\lambda^m)$, and following the same line of argument as Theorem \eqref{thm:rate:nestt-g} we eventually conclude that 
\begin{align}
\mathbb{E}[H(w^m)]\leq \frac{C(\alpha) \mathbb{E}[L^1-L^{R+1}]}{R}.
\end{align}
The proof is complete. 
\QED

\noindent\subsection{Proof of Theorem \ref{thm:lin:conv:nestte}}
Following similar line of proof as Theorem \ref{thm:lin:conv:nesttg}, we conclude that there exists some {\it finite} $\bar{r}>0$ such that for all $r>\bar{r}$ and $\bar{v}\in \mathbb{R}$ such that
\begin{align}
f(\bar{z}^r) = \bar{v}, \quad \forall~r \ge \bar{r}
\end{align}
where $\bar{z}^r=\arg\min_{z\in Z^*}\|z^r-z\|.$ Therefore, combining the fact that $\|x^{r+1}-x^r\|\to 0$, $\|z^{r+1}-z^r\|\to 0$, $\|x^{r+1}_i - z^{r+1}\|\to 0$ and $\|\lambda^{r+1}-\lambda^r\|\to 0$ [cf. \eqref{eq:difference:z:0}, \eqref{eq:difference:lambda:0}], it is easy to see that 
\begin{align}
L(\bar{z}^r, \bar{x}^r, \bar{\lambda}^r) = f(\bar{z}^r) =\bar{v}, \quad \forall~r \ge \bar{r}
\end{align}
where $\bar{x}^r, \bar{\lambda}^r$ are defined similarly as $\bar{z}^r$.

In what follows we will show that $\Delta:=L^{r+1}-\bar{v}$ decreases Q-linearly.

From the optimality condition for $z$ subproblem \eqref{eq:optcond:z}, we know that:
\begin{align}
z^{r+1}=\prox_h\left[z^{r+1}-\left(\sum_{i=1}^N -\lambda_i^r-\eta_i(z^{r+1}-x_i^r)\right) - \nabla g_0(z^{r+1})\right].
\end{align}
Therefore, $\|\tilde{\nabla}f(z^{r+1})\|$ can be bounded in the following way:
\begin{align}
\|\tilde{\nabla}f(z^{r+1})\|&=\bigg\|\prox_h\left[z^{r+1}-\left(\sum_{i=1}^N -\lambda_i^r-\eta_i(z^{r+1}-x_i^r)\right)-\nabla g_0(z^{r+1})\right] \nonumber\\
&\quad - \prox_h\left(z^{r+1}-\sum_{i=1}^N\nabla g_i(z^{r+1})- \nabla g_0(z^{r+1})\right)\bigg\| \nonumber\\
&\stackrel{\rm(a)}\leq \left\|-\sum_{i=1}^N\nabla g_i(z^{r+1})+\left(\sum_{i=1}^N -\lambda_i^r-\eta_i(z^{r+1}-x_i^r)\right)\right\|\nonumber\\
&\stackrel{\rm(b)}=\left\|\sum_{i=1}^N\left[-\nabla g_i(z^{r+1})+ \nabla g_i(x_i^{r})-\eta_i(z^{r+1}-x_i^r)\right]\right\|\nonumber\\
&\leq \sum_{i=1}^N(L_i/N+\eta_i)\|z^{r+1}-x_i^{r}\|\nonumber\\
&\leq \sum_{i=1}^N(L_i/N+\eta_i)\left(\|z^{r+1}-\hat{x}_i^{r+1}\| +\|\hat{x}_i^{r+1}-x_i^r\|\right)\nonumber\\
&= \sum_{i=1}^N(L_i/N+\eta_i)\left(\frac{1}{\alpha_i\eta_i}\|\hat{\lambda}_i^{r+1}-\lambda_i^r\| +\|\hat{x}_i^{r+1}-x_i^r\|\right)\nonumber\\
&\stackrel{\rm(c)}\leq \sum_{i=1}^N (L_i/N+\eta_i)\left(\frac{L_i}{N\alpha_i\eta_i}+1\right)\|\hat{x}_i^{r+1}-x_i^r\|
\end{align}
where in $\rm(a)$ we have used the fact that poximity operator is nonexpansive, and  in $\rm(b)$  we plugged in equation \eqref{eq:y_expression1}; in ${\rm(c)}$ we have used Lemma \ref{lemma:y1}.
Therefore, there exists $C>0$ such that we can bound the $\|\tilde{\nabla}f(z^{r+1})\|^2$ as the  following:
\begin{align}
\|\tilde{\nabla}f(z^{r+1})\|^2 \leq C\left(\|z^{r+1}-z^r\|^2 + \sum_{i=1}^N\|\hat{x}_i^{r+1}-x_i^r\|^2\right)\label{eq:bd:res}.
\end{align}
In what follows we bound $L(z^{r+1},x^r,\lambda^r)-\bar{v}$. Assume that $r>\bar{r}$, therefore  $\bar{z}^{r+1}=\arg\min_{z\in Z^*} \|z^{r+1}-z\|$  satisfies $f(\bar{z}) = \bar{v}$.  We have the following 
\begin{align}
L(z^{r+1},x^r,\lambda^r)-\bar{v} &= \sum_{i=1}^{N} g_i(x^{r}_i)+g_0(z^{r+1})+p(z^{r+1})+\sum_{i=1}^{N}\langle \lambda^{r}_i, x^{r}_i-z^{r+1}\rangle+\sum_{i=1}^{N}\frac{\eta_i}{2}\|x^{r}_i-z^{r+1}\|^2\nonumber\\
& - \sum_{i=1}^{N}  g_i(\bar{z}^{r+1})-p(\bar{z}^{r+1})-g_0(\bar{z}^{r+1})\label{eq:L:bd}.
\end{align}
From the optimality condition of the $z$ subproblem we know that:
\begin{align}
&p(z^{r+1})+ g_0(z^{r+1})+\sum_{i=1}^{N}\langle \lambda^{r}_i, x^{r}_i-z^{r+1}\rangle+\sum_{i=1}^{N}\frac{\eta_i}{2}\|x^{r}_i-z^{r+1}\|^2\nonumber\\
&\leq p(\bar{z}^{r+1})+g_0(\bar{z}^{r+1})+\sum_{i=1}^{N}\langle \lambda^{r}_i, x^{r}_i-\bar{z}^{r+1}\rangle+\sum_{i=1}^{N}\frac{\eta_i}{2}\|x^{r}_i-\bar{z}^{r+1}\|^2.
\end{align}
Rearranging the terms in the previous equation we have:
\begin{align}
&p(z^{r+1}) + g_0(z^{r+1})-p(\bar{z}^{r+1})-g_0(\bar{z}^{r+1})\nonumber\\
&\leq\sum_{i=1}^{N}\left[\langle \lambda^{r}_i, z^{r+1}-\bar{z}^{r+1}\rangle+\frac{\eta_i}{2}\|x_i^{r}-\bar{z}^{r+1}\|^2-\frac{\eta_i}{2}\|x^{r}_i-z^{r+1}\|^2\right]\label{eq:h:bd}.
\end{align}
Plugging in \eqref{eq:h:bd} in \eqref{eq:L:bd} yields the following:
\begin{align}
L(z^{r+1},x^r,\lambda^r)-\bar{v} &\leq \sum_{i=1}^{N}\left[  \frac{1}{N}(g_i(x^{r}_i)-g_i(\bar{z}^{r+1}))+\frac{\eta_i}{2}\|x^{r}_i-\bar{z}^{r+1}\|^2 + \langle \lambda^{r}_i, x_i^{r}-\bar{z}^{r+1}\rangle\right]\nonumber\\
&\stackrel{\rm(a)}=\sum_{i=1}^{N}\left[ \frac{1}{N}\langle  \nabla g_i(\tilde{x}_i), x_i^{r}-\bar{z}^{r+1}\rangle+\frac{\eta_i}{2}\|x^{r}_i-\bar{z}^{r+1}\|^2 + \langle \lambda^{r}_i, x_i^{r}-\bar{z}^{r+1}\rangle\right]\nonumber\\
&\stackrel{\rm(b)}=\sum_{i=1}^{N}\left[ \frac{1}{N}\langle  \nabla g_i(\tilde{x}_i) - \nabla g_i(x_i^{r}), x_i^{r}-\bar{z}^{r+1}\rangle+\frac{\eta_i}{2}\|x^{r}_i-\bar{z}^{r+1}\|^2\right] \nonumber\\
&\stackrel{\rm(c)}\leq \sum_{i=1}^{N}\left(L_i/N^2+\frac{\eta_i}{2}\right)\|x_i^{r}-\bar{z}^{r+1}\|^2\nonumber\\
&\leq \sum_{i=1}^{N}3(L_i/N+\eta_i)\left[\|x_i^{r}-\hat{x}^{r+1}_i\|^2+\|\hat{x}^{r+1}_i-z^{r+1}\|^2+\|z^{r+1}-\bar{z}^{r+1}\|^2\right]\nonumber\\
&\leq \sum_{i=1}^{N}3(L_i/N^2+\eta_i)\left[\left(1+\frac{L_i^2}{\alpha_i^2\eta_i^2N^2}\right)\|x_i^{r}-\hat{x}_i^{r+1}\|^2+\|z^{r+1}-\bar{z}^{r+1}\|^2\right].
\end{align}
In the above series of inequalities,  $\rm(a)$ is the result of applying mean value theorem on function $g_i$, where $\tilde{x}_i$ is a point lies in between the line segment $(x^r_i, \bar{z}^{r+1})$; In $\rm(b)$ we have used \eqref{eq:y_expression1}; $\rm(c)$ is true because of assumption A-(c) and  the fact that $\|\tilde{x_i}-\bar{z}\|\leq \|x_i^{r}-\bar{z}^{r+1}\|$; the last inequality comes from Lemma \ref{lemma:y1}. 

Overall, there exists $\sigma_2>0$ such that 
\begin{align}
L(z^{r+1},x^r,\lambda^r)-\bar{v} \leq \sigma_2\left(\sum_{i=1}^N\|\hat{x}_i^{r+1}-x_i^{r}\|^2+\|z^{r+1}-\bar{z}^{r+1}\|^2\right)\quad \forall~r\ge \bar{r}. \label{eq:L:bd2}
\end{align}
From Lemma \ref{lm:eb}, we know that
\begin{align}
\dist(z^{r+1},Z^*) = \|z^{r+1} -\bar{z}^{r+1}\|\le \tau \|\tilde\nabla f(z^{r+1})\|\label{eq:eb2}.
\end{align} 
As a result we have:
$$\Delta^{r+1}= L^{r+1}-\bar{v}\le L(z^{r+1},x^r,\lambda^r)-\bar{v},$$ 
where the inequality comes from the fact that $L^{r+1}- L(z^{r+1},x^r,\lambda^r)\le 0$ [by \eqref{eq:key:2} and the fact that $c_i<0$ for all $i$]. 
Considering the above equations we obtain:
\begin{align}
\Delta^{r+1}&\leq  \sigma_2\left(\sum_{i=1}^N\|\hat{x}_i^{r+1}-x_i^{r}\|^2+\|z^{r+1}-\bar{z}^{r+1}\|^2\right)\nonumber\\
&\stackrel{\eqref{eq:eb2}}\leq \sigma_2\sum_{i=1}^N\|\hat{x}_i^{r+1}-x_i^{r}\|^2+\sigma_2\tau\|\tilde\nabla f(z^{r+1})\|^2\nonumber\\
&\stackrel{\eqref{eq:bd:res}}\leq (\sigma_2+\sigma_2\tau C)\sum_{i=1}^{N}\|\hat{x}_i^{r+1}-x_i^{r}\|^2+\sigma_2\tau C \|z^{r+1}-z^r\|^2, \quad \forall~r\ge \bar{r}\label{eq:delta:diff}
\end{align}
From \eqref{eq:estimate:L} and \eqref{eq:delta:diff} we can further construct a $\sigma_3>0$ such that $\mathbb{E}[\Delta^{r+1}]\leq \sigma_3\mathbb{E}[L^r-L^{r+1}]$, which further implies  the following relationship:
\begin{align}
E[\Delta^{r+1}]\leq \frac{\sigma_3}{1+\sigma_3}E[\Delta^r]\label{eq:lin:conv}, \quad \forall~r\ge \bar{r}.
\end{align}
Let us set $\rho=\frac{\sigma_3}{1+\sigma_3}<1$. Thus concluding the proof. \QED

\subsection{ {Proof of Proposition \ref{prop:omega:close}}}
Applying the optimality condition on $z$ subproblem in \eqref{eq:z:prox} we have:
\begin{align}
z^{r+1}=\argmin_z h(z) + g_0(z)+\frac{\beta}{2}\|z-u^{r+1}\|^2 \label{eq:prox:omega}
\end{align}
where the variable $u^{r+1}$ is given by (cf. \eqref{eq:u})
\begin{align}
u^{r+1}=\beta\sum_{i=1}^{N}( \lambda_{i}^{r}+\eta_i x_{i}^{r+1})\label{eq:w:close:nesst}.
\end{align}
Now from  one of the key properties of NESTT-G [cf. Section \ref{sec:key:nesttg}, equation \eqref{eq:comp}],  we have that  
\begin{align}
u^{r+1} = z^r - {\beta} \left(\frac{1}{{N\alpha_{i_{r}}}}\left(\nabla g_{i_{r}} (z^r)-\nabla g_{i_{r}}(y^{r-1}_{i_{r}})\right) +\frac{1}{N} {\sum_{i=1}^{N}{\nabla g_i (y^{r-1}_i)}}{N}\right).
\end{align}

This verifies the claim. 
\QED

			\newpage
			
{
\bibliographystyle{ieee}
\bibliography{ref,biblio,distributed_opt}
}
\end{document}